\documentclass[11pt]{article}
\usepackage{amssymb, amsmath, latexsym, proof}

\title{Periodic frames}

\usepackage{authblk}

\author[1,2]{Lev D. Beklemishev\footnote{The corresponding author,  e-mail: \texttt{bekl@mi-ras.ru}.}}
\author[3]{Yunsong Wang}

\affil[1]{Steklov Mathematical Institute of Russian Academy of Sciences}
\affil[2]{National Research University Higher School of Economics, Moscow}
\affil[3]{Peking University, Beiging, China}

\newtheorem{theorem}{Theorem}
\newtheorem{lemma}{Lemma}[section]

\newtheorem{corollary}[lemma]{Corollary}
\newtheorem{proposition}[lemma]{Proposition}

\newtheorem{definition}{Definition}
\newtheorem{example}[lemma]{Example}
\newtheorem{remark}[lemma]{Remark}

\newcommand{\bl}{\begin{lemma}}
\newcommand{\el}{\end{lemma}}
\newcommand{\bt}{\begin{theorem}}
\newcommand{\et}{\end{theorem}}
\newcommand{\bcor}{\begin{corollary}}
\newcommand{\ecor}{\end{corollary}}
\newcommand{\bp}{\proof{.}}
\newcommand{\ep}{\eop}
\newcommand{\bpr}{\begin{proposition}}
\newcommand{\epr}{\end{proposition}}
\newcommand{\brem}{\begin{remark} \em}
\newcommand{\erem}{\end{remark}}
\newcommand{\bd}{\begin{definition} \em}
\newcommand{\ed}{\end{definition}}
\newcommand{\bex}{\begin{example} \em
}
\newcommand{\eex}{\end{example}}
\newcommand{\beq}{\begin{equation} }
\newcommand{\eeq}{\end{equation}}

\newcommand{\bi}{\begin{itemize}
  }
\newcommand{\ei}{\end{itemize}}
\newcommand{\ben}{\begin{enumerate} }
\newcommand{\een}{\end{enumerate} }

\newenvironment{enumr}{

\begin{enumerate}     }{\end{enumerate}

}
\newcommand{\benr}{\begin{enumr}
  }
\newcommand{\eenr}{
\end{enumr}}

\newcommand{\bs}{\bigskip}
\newcommand{\ignore}[1]{}

\newcommand{\al}[1]{\forall #1\:}

\newlength{\hilflh}

\renewcommand{\emptyset}{\varnothing}

\newcommand{\cL}{{\mathcal L}}
\newcommand{\cB}{{\mathcal B}}

\newcommand{\cT}{{\mathcal T}}
\newcommand{\cP}{{\mathcal P}}

\newcommand{\cA}{{\mathcal A}}
\newcommand{\cM}{{\mathcal M}}
\newcommand{\cV}{{\mathcal V}}

\newcommand{\cW}{{\mathcal W}}

\newcommand{\cX}{{\mathcal X}}

\newcommand{\cS}{\mathcal{S}}

\newcommand{\ga}{\alpha}
\newcommand{\gb}{\beta}
\newcommand{\gd}{\delta}
\renewcommand{\ge}{\varepsilon}
\newcommand{\gl}{\lambda}

\newcommand{\gs}{\sigma}
\newcommand{\gy}{\gamma}
\newcommand{\gw}{\omega}

\renewcommand{\phi}{\varphi}

\newcommand{\imp}{\rightarrow}

\newcommand{\GL}{\mathbf{GL}}
\newcommand{\GLP}{\mathbf{GLP}}

\newcommand{\la}{\langle}
\newcommand{\ra}{\rangle}

\renewcommand{\models}{\vDash}      %amsmodels

\newcommand{\nmodels}{\nvDash}

\newcommand{\On}{\mathrm{On}}

\newcommand{\Lim}{\mathrm{Lim}}

\newcommand{\J}{\mathbf{J}}

\renewcommand{\leq}{\leqslant}
\renewcommand{\geq}{\geqslant}

\newcommand{\eop}{$\Box$ \protect\par \addvspace{\topsep}}
\newcommand{\proof}[1]{\protect\par\addvspace{\topsep}\noindent {\bf Proof#1}}
\begin{document}

\maketitle

\begin{abstract}
Polymodal provability logic GLP is incomplete w.r.t.\ Kripke frames. It is known to be complete w.r.t.\  topological semantics, where the diamond modalities correspond to topological derivative operations. However, the topologies needed for the completeness proof are highly non-constructive.
The question of completeness of GLP w.r.t.\ natural scattered topologies on ordinals is dependent on large cardinal axioms of set theory and is still open. So far, we are lacking a useable class of models for which GLP is complete.

In this paper we define a natural class of countable general topological frames on ordinals for which GLP is sound and complete. The associated topologies are the same as the ordinal topologies introduced by Thomas Icard. However, the key point is to specify a suitable algebra of subsets of an ordinal closed under the boolean and topological derivative operations. The algebras we define are based on the notion of a periodic set of ordinals generalizing that of an ultimately periodic binary omega-word.
\end{abstract}

\newcommand{\VV}{\mathcal{V}}
\newcommand{\ml}[1]{\mathrm{#1}}
\newcommand{\HH}{\mathcal{H}}
\newcommand{\D}{d}
\newcommand{\tto}{\twoheadrightarrow}
\newcommand{\Log}{\mathrm{Log}}
\renewcommand{\GL}{\mathrm{GL}}
\renewcommand{\GLP}{\mathrm{GLP}}
\renewcommand{\J}{\mathrm{J}}

\newcommand{\lr}[1]{(#1)}

\section{Introduction}
Polymodal provability logic $\GLP$ was introduced by Giorgi  Japaridze~\cite{Dzh86,Boo93} as the logic of the extensions of Peano arithmetic by $n$-fold applications of the omega-rule, for each $n<\gw$. Later, Konstantin Ignatiev extended Japaridze's work and showed that $\GLP$ is sound and complete w.r.t.\ a wider class of arithmetical interpretations by any sufficiently strong sequence of sound provability predicates~\cite{Ign93}. He also obtained normal forms and introduced a universal Kripke model for the variable-free fragment of $\GLP$ that played a prominent role in later studies.

More recently, $\GLP$ has found interesting applications in
proof-theoretic analysis of arithmetic \cite{Bek04,Bek05,Bek06,Joo15a,JF18,PW21,JooFPP22,BekPakh22}
which stimulated further interest in the study of modal-logical
properties of $\GLP$ \cite{BJV,Shap08,Ica09,Bek10,JF13,Pakh14,Sham11,Sham15}.

\ignore{
For such
applications, the algebraic language appears to be more natural and a different choice of the interpretation of the provability
predicates is needed. The relevant structures have been introduced
in \cite{Bek04} under the name \emph{graded provability algebras}.
}

One obstacle in the study of $\GLP$ is its lack of Kripke completeness. It is easy to see that in any Kripke frame validating the axioms of $\GLP$ the relations corresponding to all modalities but $[0]$ have to be empty. Thereby, Kripke completeness of $\GLP$ fails in a very strong way. Over the years there has been a steady effort to develop some natural and usable models for this logic.

A sufficiently rich class of \emph{Kripke models} (rather than frames), for which $\GLP$ is sound and complete, has been introduced in~\cite{Bek10}. The description of this class, although constructive, is rather complicated. In particular, appropriate models are obtained from finite ones by a double limit construction.
Thereafter, attention of researchers was drawn to \emph{neighborhood} or \emph{topological models} of $\GLP$.

Under the standard topological interpretation of provability logic $\GL$, propositional variables are interpreted as subsets of a (scattered) topological space, boolean connectives correspond to boolean operations, whereas the modality $\Diamond$ is interpreted as the topological derivative operation mapping a set $A$ to the set of its limit points. This interpretation originated in the work of H.~Simmons~\cite{Sim75} and L.~Esakia~\cite{Esa81} who isolated the role of scatteredness, the topological counterpart of L\"ob's axiom. Later, Abashidze~\cite{Aba85} and Blass~\cite{Bla90} proved the completeness of $\GL$ for the topological interpretation w.r.t.\ the interval topology on any ordinal $\Omega$ provided $\Omega\geq\gw^\gw$. This result served as a basis for various generalizations. We refer to \cite{BekGab14} for a comprehensive survey of topological models of provability logic.

Topological interpretation of $\GLP$ turned out to be a rich area of study on its own with a number of significant contributions.
Thomas Icard~\cite{Ica09} suggested natural topological models for the variable-free fragment of $\GLP$ based on Ignatiev's work. These models fail to characterize full $\GLP$, however, they are somewhat simpler than Ignatiev's Kripke frame. Icard's topologies on ordinals will play a major role in this paper.

In \cite{Bek09,BBI09} the notion of \emph{$\GLP$-space}, that is, a sound topological model of $\GLP$, was introduced and it was shown that the bimodal fragment of $\GLP$ is topologically complete.
In \cite{Bek10b}, this has been improved by showing that the bimodal fragment of $\GLP$ is complete w.r.t.\ a natural pair of ordinal topologies (assuming the set-theoretic hypotheses $V=L$ or Jensen's square principles). This paper was based on an earlier result by Andreas Blass~\cite{Bla90} who showed that $\GL$ is sound and complete w.r.t.\ the club filter topology on $\aleph_\gw$.

The extension of this result to more than two modalities turned out to be a significant problem. The ordinals bear a natural sequence of topologies making it a $\GLP$-space, the so-called \emph{ordinal $\GLP$-space}. This model and the underlying topologies were introduced in~\cite{Bek09,BBI09}. However, whether $\GLP$ is complete for this model is still an open question, even under set-theoretic assumptions such as $V=L$ and large cardinal axioms. It was noticed that in order to guarantee that all the topologies in that sequence are at least non-discrete --- a necessary condition for completeness --- one really needs large cardinal assumptions. A natural sufficient condition for the non-discreteness is that $\Pi_n^1$-indescribable cardinals exist, for each $n<\gw$.
More exactly, the limit points of the topologies of the ordinal GLP-space correspond to the so-called \emph{two-fold $d_n$-reflecting ordinals}~\cite{Bek09,BekGab14} or equivalently, using the terminology of~\cite{Bag12,BaMaSa12}, to \emph{$n$-simultaneusly-stationary reflecting cardinals}, for $n<\gw$.

The question of consistency strength of $n$-(simultaneusly)-stationary reflecting cardinals and its generalization to transfinite iterations turned out to be fruitful from a set-theoretic point of view and stimulated a number of papers starting from the work of Joan Bagaria~\cite{Bag12}, see~\cite{BaMaSa12,BaMaMa20,Ba19}. See also~\cite{BekGab14} for a survey of set-theoretic aspects of $\GLP$.

Despite these set-theoretic complications the topological completeness of $\GLP$ was established in~\cite{BekGab13} on the basis of ZFC by a substantially new construction. However, the constructed topological spaces were not natural and, in fact, were obtained with the heavy use of the axiom of choice. The technical difficulties involved were also significant. Later this result has been extended and improved in various ways. In particular, Fern\'andez-Duque~\cite{Fern14-top} extended it to a version of $\GLP$ in the language with countably many modalities. For the uncountably many modalities the problem is still open. Working towards its solution, Aguilera~\cite{Agu23} proved the topological completeness of $\GLP$ w.r.t.\ a model of generalized Icard topologies over any sufficiently tall scattered space. Shamkanov~\cite{Sham20} found an alternative proof of topological completeness of $\GLP$ avoiding the intermediate use of Kripke models and proceeding directly from the algebraic models. He also established the \emph{strong completeness} of $\GLP$ w.r.t.\ topological semantics generalizing an earlier result of Aguilera and Fern\'andez-Duque~\cite{AguFer17} for the case of G\"odel--L\"ob logic.

In this paper we consider a more general kind of topological semantics for $\GLP$ called \emph{general topological (derivative) frames}. General topological frames for a modal logic relate to its topological semantics in the same way as general Kripke frames relate to Kripke semantics. These are polytopological spaces $X$ equipped with an algebra of admissible subsets $\cA\subseteq \cP(X)$. The algebra is supposed to be closed under the boolean operations and the derivative operations corresponding to the topologies of the structure. Propositional formulas are then interpreted as elements of $\cA$. %This makes general topological frames essentially the same thing as subalgebras of the $\GLP$-algebras represented as the algebras of subsets of a set $X$.

Such structures have been considered in the modal logic literature before, e.g., the authors of~\cite{BezGabLuc15} consider general topological frames for S4 based on topological closure rather than derivative operators. General topological derivative frames have been introduced by David Fern\'andez-Duque~\cite{Fern14-top} under the name \emph{ambiences}. According to~\cite{Fern14-top}, a \emph{provability ambience} is a general topological frame satisfying the axioms of $\GLP$ in which the sequence of topologies is increasing and the topologies are scattered.\footnote{In principle, a general topological frame satisfying the axioms of $\GLP$ need not be topologically scattered. However, in this paper as well as in \cite{Fern14-top} we only deal with scattered topologies, that is, de facto ambiences.} We will call provability ambiences \emph{$\GLP$-ambiences}.

\ignore{
In~\cite{BBI09} a related notion of \emph{descriptive $\GLP$-frame} was introduced and an (easy) Stone-type representation theorem for $\GLP$-algebras w.r.t.~descriptive $\GLP$-frames was proved (cf.~Lemma 2.11 therein). Descriptive $\GLP$-frame is a particular kind of general topological frame, where the algebra of admissible sets is the algebra of clopens of a Stone space $X$, and the topologies correspond to upsets of some binary relations on $X$. In general, such a frame need not be topologically scattered.}

In this paper we define explicit and natural general $\GLP$-ambiences for which the logic $\GLP$ is complete. We exhibit structures of this kind in the polytopological space of the ordinal $\ge_0$, equipped with the sequence of Icard's topologies. (In \cite{Fern14-top} such ambiences are called \emph{shifted Icard ambiences}.) In our treatment the algebra of admissible subsets on $\ge_0$ is defined by adopting a very natural idea of periodicity.

Periodic transfinite binary words, or equivalently periodic subsets of the ordinals, have appeared in the literature before. On the one hand, they are considered in the context of the combinatorics of words, as a natural generalization of periodic $\gw$-words (see, for example,~\cite{CaCh01}).
On the other hand, they appeared in the context of the studies of monadic second-order theories of well-orders, starting from the work of J.R.~B\"uchi~\cite{Buc65,Buc73}. In fact, MSO definable sets on $\gw^\gw$ happen to be precisely the  \emph{hereditarily periodic sets} of this paper. (It seems to be the most natural notion of a periodic set of ordinals that one can come up with.) However, these sets are not enough for the completeness of $\GLP$ in the language with more than one modality. So, we will generalize them to what we call \emph{$n$-periodic sets}, for each $n<\gw$, and to arbitrary ordinals.

The main result of this paper is the completeness theorem for $\GLP$ w.r.t.\ the ordinal $\ge_0$ (or any larger ordinal) equipped with Icard's topologies and the algebra of $n$-periodic sets. This structure gives the best solution that we have at the moment to the problem of finding an explicit and usable complete semantics for the provability logic $\GLP$.

The paper is organized as follows. In Section 2 we present necessary background results on $\GLP$, its algebraic and neighborhood semantics, and general topological frames. In Section 3 we introduce periodic transfinite words and sets of ordinals. We prove that periodic sets form a Magari algebra. In Section 4 we introduce $n$-periodic sets and the corresponding general topological frames. In Section 5 we borrow some techniques from~\cite{BekGab13} reducing the question of completeness of $\GLP$ w.r.t.\ periodic frames to a covering construction. Some necessary intermediate lemmas from~\cite{BekGab13} are presented in the appendix to this paper. Finally, in Section 6 we prove the covering lemma and thereby complete the proof of our main result.

The authors started working on this topic together in the context a Master's thesis project of Yunsong Wang at the ILLC of the University of Amsterdam. In his Thesis~\cite{YW23} the notions of hereditarily periodic set and of 1-periodic set were defined, and the corresponding completeness theorem for the bimodal fragment of $\GLP$ was proved. After the completion of this work, the first author of the paper reworked the whole setup and extended the result to the full language of $\GLP$. The current paper is written up by Lev Beklemishev.

We are grateful to Nick Bezhaishvili, a co-supervisor of the Master's Thesis of Yunsong Wang, and to David Fern\'andez-Duque and Guram Bezhanishvili for valuable comments.

\section{Preliminaries}
Here we introduce some of the background notions for this paper.

\subsection{Provability logic $\GLP$}

$\GLP$ is a propositional modal logic formulated in the language with infinitely many
modalities $[n]$, for all $n\in\gw$. As usual, we abbreviate $\neg[n]\neg \phi$ by $\la n\ra\phi$. $\GLP$ is
given by the following axiom schemata and rules:
\begin{description}
\item[Axioms:]
\begin{enumr}
\item Boolean tautologies; \item $[n](\phi\imp\psi)\imp
([n]\phi\imp [n]\psi)$; \item $[n]([n]\phi\imp \phi)\imp [n]\phi$;
\item $\la m\ra\phi
\imp [n]\la m\ra\phi$, for $m<n$.
\item $[m]\phi\imp [n]\phi$, for
$m\leq n$; \label{mon}
\end{enumr}
\item[Rules:] modus ponens, $\phi\vdash [n]\phi$.
\end{description}

The first three axioms are the well-known axioms of the G\"odel--L\"ob provability logic $\GL$ (stated for modality $[n]$). We identify $\GL$ with the fragment of $\GLP$ restricted to the language with modality $[0]$. 

\ignore{
While $\GLP$ is Kripke incomplete, the following subsystem $\J$ turned out to be useful in the study of its semantics. An axiomatization of $\J$ is obtained by replacing Axiom (\ref{mon}) by two theorems of $\GLP$:

\benr
\item $ [m]\phi
\imp [n][m]\phi$, for $m\leq n$;
\item $[m]\phi\to [m][n]\phi$, for $m\leq n$.
\eenr

As shown in \cite{Bek10}, $\J$ is complete w.r.t.\ a suitable class of finite Kripke frames.
}

\subsection{GLP-algebras}
Algebraic models of $\GLP$ are called $\GLP$-algebras. These are boolean algebras with operators $\la n\ra$, for all $n\in\gw$, satisfying the following modal identities:

\benr
\item $\la n\ra(x\lor y)=\la n\ra x \lor \la n\ra y$;
\quad $\neg\la n\ra\bot$;
\item $\la n\ra x = \la n\ra(x\land \neg\la n\ra x)$;
\item $\la n\ra x\land \la m\ra y = \la n\ra (x\land \la m\ra y)$, for $m<n$;
\item $\la n \ra x \land \la m\ra x = \la n\ra x$, for $m\leq n$.
\eenr

Identities (i) and (ii), for any modality $\la n\ra$, define the variety of \emph{Magari algebras}, which are the algebraic models of $\GL$. The last two identities correspond to the remaining axioms of $\GLP$.

\ignore{
Let $t_\phi$ denote the $\GLP$-algebra term corresponding to a modal formula $\phi$. The \emph{logic} of a $\GLP$-algebra $\cM$ is the set
$$\Log(\cM):=\{\phi:\cM\models \al{\vec x} t_\phi(\vec x)=\top\}.$$
We say that a logic $L$ is complete for a class of algebras $\mathbf{C}$, if $L=\bigcap_{\cM\in\mathbf{C}}\Log(\cM)$.}

A \emph{valuation} on a $\GLP$-algebra $\cM$ is a map $v$ from propositional variables to $\cM$. It is inductively extended to the set of all $\GLP$-formulas by the rules: $v(\phi\land \psi)=v(\phi)\land v(\psi)$, $v(\neg\phi)=\neg v(\phi)$, and $v([n]\phi)=\neg\la n\ra\neg v(\phi)$.

 The \emph{logic} of $\cM$ is the set
$\Log(\cM)$ of all formulas $\phi$ such that $v(\phi)=\top$, for all valuations $v$ in $\cM$.
We say that a logic $L$ is \emph{complete} for a class of algebras $\mathbf{C}$, if $L=\bigcap_{\cM\in\mathbf{C}}\Log(\cM)$.

On general grounds, the logic $\GLP$ is complete for the class of all $\GLP$-algebras, and $\GL$ is complete for the class of all Magari algebras.

\ignore{
\bpr
$\GLP\vdash \phi(\vec p)$ iff $\cM\models \al{\vec x} T_\phi(\vec x)=\top$, for all $\GLP$-algebras $\cM$.
\epr
}

\subsection{Neighborhood frames and $\GLP$-spaces}

Let $L$ be a normal modal logic in the language with modalities $[i]$, for $i\in I$. We call a \emph{neighborhood frame for $L$} a structure $(X,(\gd_i)_{i\in I})$ where each $\gd_i$ is an operator $\gd_i:\cP(X)\to \cP(X)$ and the algebra $(\cP(X), (\gd_i)_{i\in I})$ satisfies the axioms of $L$ under the interpretation of boolean connectives as the usual set-theoretic operations and $\la i\ra$ as $\gd_i$.

Esakia~\cite{Esa81} has shown that in the case of G\"odel--L\"ob logic $\GL$ every such $\gd$ coincides with the topological derivative operation of a (uniquely defined) scattered topology on $X$. Recall that a topological space $(X,\tau)$ is \emph{scattered} if each nonempty subspace of $X$ has an isolated point. Topological derivative operation $d_\tau$ maps a subset $A\subseteq X$ to its set of limit points $d_\tau(A)$. Esakia's result justifies us calling neighborhood frames for $\GL$ \emph{topological (derivative) frames}.
\ignore{In this paper we will only deal with the derivative interpretation of modalities.}

The same applies to neighborhood frames for $\GLP$ which correspond to \emph{GLP-spaces}~\cite{BBI09}.

\bd
$(X,(\tau_i)_{i\in \omega})$ is a \emph{$\GLP$-space} if, for each $n\in\omega$:
        \begin{itemize}
            \item[(i)] $(X,\tau_n)$ is a scattered topological space;
            \item[(ii)] $\tau_n\subseteq\tau_{n+1}$;
            \item[(iii)] For each $A\subseteq X$, $d_{\tau_n}(A)$ is $\tau_{n+1}$-open.
        \end{itemize}
\ed
It is easy to see that $(X,(\tau_i)_{i\in \omega})$ is a $\GLP$-space iff $(X,(d_{\tau_i})_{i\in \omega})$ is a neighborhood frame for $\GLP$ iff $(\cP(X),(d_{\tau_i})_{i\in \omega})$ is a $\GLP$-algebra.

\subsection{General topological frames and ambiences}
Ambiences have been introduced by Fern\'andez-Duque in~\cite{Fern14-top} as a generalization of $\GLP$-spaces. The same notion is called \emph{general topological frame} in~\cite{YW23}.
\bd
     An \emph{ambience} is a tuple $\lr{X,(\tau_i)_{i\in I},\mathcal{A}}$ where $\lr{X,(\tau_i)_{i\in I}}$ is a poly-topological space and $\mathcal{A}\subseteq\mathcal{P}(X)$ is closed under the boolean operations and $d_{\tau_i}$, for each $i\in I$. We call $\cal A$ the algebra of \emph{admissible sets} of the ambience.
\ed
Similarly to neighborhood frames, with each ambience $\mathcal{X}=\lr{X,(\tau_i)_{i\in I},\mathcal{A}}$ we associate the boolean algebra with operators $(\cA,(d_{\tau_i})_{i\in I})$.
Thus, an ambience is essentially the same thing as a subalgebra of the algebra associated with a topological derivative frame.  By definition, the \emph{logic of $\mathcal X$}, denoted $\Log(\mathcal{X})$, is the logic of this algebra.

\bd
        An ambience $\lr{X,(\tau_n)_{n\in \gw},\mathcal{A}}$ is a \emph{$\GLP$-ambience} (or \emph{provability ambience}~\cite{Fern14-top}) if, for each $n\in\gw$:
        \begin{itemize}
            \item[(i)] $\lr{X,\tau_n}$ is a scattered topological space;
            \item[(ii)] $\tau_n\subseteq\tau_{n+1}$;
            \item[(iii)] For each $U\in {\cal A}$, $d_{\tau_n}(U)$ is $\tau_{n+1}$-open.
        \end{itemize}
\ed

If $\mathcal X$ is a $\GLP$-ambience then its algebra is a $\GLP$-algebra and $\Log(\mathcal{X})$ contains $\GLP$.  $\GLP$-spaces are $\GLP$-ambiences in which $\cA$ is the algebra of all subsets of $X$.

\brem The requirements stated in items (i) and (ii) are actually stronger than is necessary for the soundness of $\GLP$. For example, L\"ob's  axiom will be satisfied if the topologies are `scattered' w.r.t.\ the admissible subsets only.
\erem

Let $(\mathcal X_j)_{j\in J}$ be a family of ambiences, $\mathcal X_j=(X_j,(\tau_{ij})_{i\in I},\cA_j)$. The \emph{sum} $\mathcal X=\bigsqcup_{j\in J} \mathcal X_j $ is defined in a natural way.  We define $\mathcal X=(X,(\tau_i)_{i\in I},\cA)$ so that $X$ is the disjoint union of sets $X_j$. A subset $U\subseteq X$ is $\tau_i$-open if $U\cap X_j\in\tau_{ij}$, for all $j\in J$, and $\cA$ consists of all sets $A\subseteq X$ such that $A\cap X_j\in \cA_j$, for all $j\in J$.

Clearly, in this case $(\cA,(d_{\tau_i})_{i\in I})$ is isomorphic to the product of the derivative algebras $(\cA_j,(d_{\tau_{ij}})_{i\in I})$. Since varieties of algebras are closed under products, if all $\mathcal X_j$ are $\GLP$-ambiences, then so is $\mathcal X$.

\bd
Let $\GLP_n$ denote the fragment of $\GLP$ in the language with modalities $[i]$ for $i\leq n$. Thus, $\GL$ is a notational variant of $\GLP_0$. The notions of $\GLP$-algebra, $\GLP$-space, $\GLP$-ambience, etc.\ have obvious analogues for $\GLP_n$, they will occasionally be used below.
\ed

\subsection{Ranks and $d$-maps}

Let $(X,\tau)$ be a scattered space. The \emph{Cantor--Bendixson sequence} of subsets of $X$ is defined by transfinite iteration of the derivative operation $d_\tau$ on $X$: $d_\tau^0 X=X$;
\quad $d_\tau^{\alpha+1} X=d_\tau(d_\tau^{\alpha} X)$ and
$d_\tau^\alpha X=\bigcap\limits_{\beta<\alpha}d_\tau^\beta X$ if $\alpha\in\Lim$.

Notice that if $X$ is scattered, all $d_\tau^\ga X$ are closed and hence decreasing. In fact, Cantor proved that
$(X,\tau)$ is scattered iff $d_\tau^\ga X = \emptyset$, for some $\ga$.

The \emph{rank function} $\rho_\tau:X\to \On$ is defined by
$$\rho_\tau(x):=\min\{\ga: x\notin d_\tau^{\ga+1}(X)\}.$$
$\rho_\tau$ maps $X$ onto an ordinal $\rho_\tau(X)$ called the Cantor--Bendixson rank of $X$.

Let $\Omega$ be an ordinal. Its \emph{left
topology} $\tau_\leftarrow$ is given by intervals of the form $[0,\ga)$, for all $\ga\in\Omega$. In this case we have $\rho(\ga)=\ga$, for all $\ga$.

The \emph{interval topology} on $\Omega$ is generated by $\{0\}$ and all intervals of the form $(\ga,\gb)$ for $\ga<\gb\leq\Omega$. Its rank function is the function $\ell$ defined by
$$\ell(0)=0; \quad \ell(\ga)=\gb \text{ if $\ga=\gy+\gw^{\gb}$, for some
$\gy$, $\gb$.}$$

Next we introduce $d$-maps that are very useful in the study of topological derivative semantics. A map $f:X\to Y$ between topological spaces is called a \emph{d-map} if $f$ is continuous, open and \emph{pointwise discrete}, that is,
$f^{-1}(y)$ is a discrete subspace of $X$ for each $y\in Y$.
$d$-maps satisfy the properties expressed in the
following proposition (see \cite{BEG05}).

\bpr\label{l:d-maps}
\begin{enumr}
\item $f^{-1}(d_Y(A))= d_X(f^{-1}(A))$, for any $A\subseteq Y$;
\item $f^{-1}:(\cP(Y),d_Y)\to (\cP(X),d_X)$ is a homomorphism of modal algebras;
\item If $f$ is onto, then $f^{-1}$ is an embedding, and hence
\item $\Log((\cP(X),d_X))\subseteq \Log((\cP(Y),d_Y))$.
\end{enumr}
\epr

\ignore{
From lemma~\ref{l:d-maps}(i) we easily obtain the following corollary
by transfinite induction.

\bcor\label{CB-pres} Suppose $f:X\to Y$ is a $d$-map. Then, for each
ordinal $\ga$, $d_X^\ga X= f^{-1}(d_Y^\ga Y)$. \ecor
}

When considering analogs of $d$-maps between ambiances, the relevant notion is \emph{algebra preserving $d$-map}.

Let  $\mathcal{X}=(X,(\tau_i)_{i\in I},\cA)$ and $\mathcal{Y}=(Y,(\gs_i)_{i\in I},\cB)$  be ambiances. A map  $f:\mathcal{X}\to \mathcal{Y}$ is called \emph{algebra preserving} if
$f^{-1}(B)\in \cA$, for each $B\in \cB$. It is called a $d$-map, if it is a $d$-map between $(X,\tau_i)$ and $(Y,\gs_i)$, for all $i\in I$.

\bpr
If $f:\mathcal{X}\to\mathcal{Y}$ is an algebra preserving $d$-map, then
\benr \item $f^{-1}:(\cB,(d_{\gs_i})_{i\in I})\to (\cA,(d_{\tau_i})_{i\in I})$ is a homomorphism;
\item If $f$ is onto, then $f^{-1}:\cB\to \cA$ is an embedding, and hence
\item
$\Log(\mathcal{X})\subseteq \Log(\mathcal{Y})$.
\eenr
\epr

The following proposition~\cite{BekGab14} states that the rank function, when the ordinals
are equipped with their left topology, becomes a $d$-map. It is also
uniquely characterized by this property.

\bpr \label{rank-dmap} Let $\Omega:=\rho_\tau(X)$ taken
with its left topology. Then
\begin{enumr}
\item $\rho_\tau:X\tto \Omega$ is an onto $d$-map;
\item If $f:X\to \gl$ is a $d$-map, where $\gl$ is an ordinal with its left topology, then $f(X)=\Omega$ and $f=\rho_\tau$.
\end{enumr}
\epr

\subsection{Icard's topologies}
\ignore{
A natural topological model for the variable-free fragment of $\GLP$ was introduced by Thomas Icard \cite{Ica09}. Considered as a poly-topological space it is not a GLP-space, thus, not a model of the full $\GLP$. However, it is sound and complete for the variable-free fragment of $\GLP$.
}

Let $\Omega$ be an ordinal. Icard's topologies $\vartheta_n$, for each $n\in\gw$, are defined as follows.\footnote{Our $\vartheta_n$ is Icard's topology $\theta_{n+1}$.}
\bd
Let $\vartheta_0$ be the interval topology, and let $\vartheta_n$ be generated as a subbase by $\vartheta_0$ and all sets of the form $$U^m_\gb(\Omega):=\{\ga\in \Omega:\ell^m(\ga)>\gb\},$$ for $m\leq n$. Here $\ell^m$ denotes the $m$-th iterate of the $\ell$ function. We let $d_n$ and $\rho_n$ denote the derivative operator and the rank function for $\vartheta_n$, respectively.
\ed

Clearly, $(\vartheta_n)_{n\in\gw}$ is an increasing sequence of scattered topologies. The following characterizations are from~\cite{BekGab14} and are easy to establish on the basis of Proposition~\ref{rank-dmap}.

\bpr \begin{enumr}
\item $\ell:(\Omega,\vartheta_{n+1})\to (\Omega,\vartheta_n)$ is a $d$-map;
\item $\vartheta_{n+1}$ is the coarsest topology $\tau$ on $\Omega$ such that $\ell:(\Omega,\tau)\to (\Omega,\vartheta_n)$ is continuous;
\item $\ell^{n+1}$ is the rank function of $\vartheta_n$, that is, $\rho_n=\ell^{n+1}$;
\item $\vartheta_{n+1}$ is generated by $\vartheta_n$ and $\{d_n^{\ga+1}(\Omega): \ga<\rho_n(\Omega)\}$.
\end{enumr}
\epr

Notice that (iv) allows one to define the analogues of Icard's topologies on any scattered space.

\ignore{
\bl\label{discr}
$\ell: \gw^\Omega\to\Omega$ is $\vartheta_0$-pointwise discrete.
\el
\bp\ We have to show that $\ell^{-1}(\gb)$ is discrete, for each $\gb$. Suppose $\ell(\ga)=\gb$ and $\ga>0$, then $\ga=\ga_0+\gw^\gb$ in Cantor normal form. We claim that the interval $(\ga_0,\ga)$ does not intersect $\ell^{-1}(\gb)$. Any ordinal in that interval has the form $\ga_0+\gy$ for some $0<\gy<\gw^\gb$. This entails that $\gy$ is a sum of $\omega$-powers smaller than $\gw^\gb$ and hence not in $\ell^{-1}(\gb)$. \ep
}

\section{Periodic sets of ordinals}

\subsection{Transfinite words and periodic sets}

Here we develop some basic theory of transfinite words. We will only consider binary words in this paper, although the notions can be easily generalized to any alphabet.

\bd
A \emph{(transfinite) word} is map $A:\alpha \to \{0,1\}$ where $\alpha\in On$. The ordinal $\ga$ is called the \textit{length} of $A$ and denoted $|A|$. We also call $A$ an \emph{$\ga$-word}. The empty word (the only word of length $0$) is denoted $\Lambda$.

\emph{Concatenation} of transfinite words $A$ and $B$ is defined in a natural way and is denoted $AB$. This operation obviously satisfies  $|AB|=|A|+|B|$, $A(BC)=(AB)C$.

A word $A$ is an \emph{initial segment} of $B$ if there is a $C$ such that $AC=B$. Similarly, $A$ is an end segment of $B$ if there is a $C$ such that $CA=B$.

\emph{Iteration} $A^\gb$ of a transfinite word $A$ is defined by transfinite induction on $\gb$: $A^0:=\Lambda$, $A^{\gb+1}:=A^\gb A$, $A^\lambda:=\bigcup_{\beta<\lambda}A^\beta$ if $\lambda\in \Lim$.

Clearly, $|A^\gb|=|A|\cdot \gb$. We also have $A^\gb A^\gy=A^{\gb+\gy}$ and $(A^\gb)^\gy=A^{\gb\gy}$.
\ed

Next we are going to introduce one of the main concepts in this paper, that of a periodic word. We do not intend to treat this notion in full generality, but we deal with a somewhat restricted concept suitable for our purposes. In particular, in our version a periodic word must always contain $\gl$-many periods, for some limit ordinal $\gl$. Thus, for example, a word of the form $AAA$ is not always periodic, but $A^\gw$ is.

\bd A word
$X$ is \emph{periodic} if $X=A^\lambda$ for some $\lambda\in \Lim$ and some $A\neq \Lambda$. In this case, $A$ is called a \emph{period} of $X$ and $\gl$ is called the \emph{exponent}.

$X$ is \emph{ultimately periodic} (u.p.) if $X=BY$ for some $Y,B$ with $Y$ periodic. Periods of $X$ are, by definition, periods of $Y$ for any such representation. $B$ is called an \emph{offset} of $X$.
\ed

We remark that these definitions agree with the standard concepts of periodic and ultimately periodic $\gw$-words.

Given an ordinal $\Omega$, there is a natural one-to-one correspondence (depending on $\Omega$) between transfinite binary words of length $\Omega$ and subsets of $\Omega$: We associate with a set $X\subseteq\Omega$ its characteristic function $\Omega\to \{0,1\}$, and with a word $A:\Omega\to \{0,1\}$ the set $\{\ga\in \Omega: A(\ga)=1\}$. Thereby, all the notions introduced so far for transfinite words are also meaningful for sets.

\bd We call a set $X$
\emph{periodic/u.p.\ in} $\Omega$ if so is the $\Omega$-word associated with the subset $X\cap \Omega$ of $\Omega$.
\ed

Thus, we extend the terminology to sets $X$ not necessarily contained in $\Omega$. Clearly, this makes sense only for limit ordinals $\Omega$ (otherwise there are no periodic sets in $\Omega$).

\subsection{Periods}
Our next goal is to show that u.p.\ subsets of $\Omega$ form a Boolean algebra. To this end, we need to look into the structure of periods of u.p.\ sets and to establish a useful \emph{synchronization lemma}.

In general, period of a word, if it exists, is not unique. For example, we obviously have
\bl \label{double}
    If $A$ is a period of $X$, then so is $A^n$, for all $n<\gw$.
\el
\bp\ $(A^n)^\gl = A^{n\gl} = A^\gl$, since $\gl\in\Lim$. \ep

An offset of an u.p.\ word is also not unique. In fact, if $B$ is an offset of $X$, then so is the extension of $B$ by any number of  periods as long as it is a proper initial segment of $X$. We can always assume the exponent of an u.p.\ word to be a power of $\omega$. Indeed, if $\gl=\gb+\gw^\ga$ with $\ga>0$ then $A^\gl=A^\gb A^{\gw^\ga}$ and one can extend the offset by $A^\gb$. We will use this principle and Lemma~\ref{double} without reference.

\bl \label{nfperiod}
    If $X$ is u.p., then $X$ has a period of length $\omega^\beta\cdot n$, for some $\beta,n$ with $0<n<\omega$.
\el

\bp\ Let $X=BA^\alpha$, $\alpha\in \Lim$. Consider two cases:

\textsc{Case 1:} $\alpha=\omega\alpha'$ with $\alpha'\in \Lim$. Then $A^\alpha=(A^\omega)^{\alpha'}$ and $|A^\omega|=|A|\cdot \omega = \omega^\beta$, for some $\gb$.

\textsc{Case 2:} $\alpha=\alpha''+\omega$. Then $A^\alpha=CA^\omega$, for some $C$. Assume $|A|=\omega^{\alpha_1}\cdot n +\gamma$ in Cantor normal form and $A=A_1A_2$ such that $|A_1|= \omega^{\alpha_1}\cdot n$ and $|A_2|= \gamma$.
Then $A^\omega = (A_1A_2)^\omega= A_1(A_2A_1)^\omega$ where $|A_2A_1|=\gamma+\omega^{\alpha_1}n=\omega^{\alpha_1}n$.
\ep

We say that $\ga$ is a (right) multiple of $\gb$ if $\ga=\gb\cdot \gy$, for some $\gy>0$.
\bl\
Multiples of $\gw^\gb$ are precisely the ordinals $\ga$ such that $\ell(\ga)\geq \gb$.
\el

\bp\ Let $\ga=\gw^{\ga_1}+\cdots +\gw^{\ga_n}$ in Cantor normal form with $\ga_1\geq \cdots \geq \ga_n$. If $\ell(\ga)=\ga_n\geq \gb$ then there are $\ga_1'\geq \cdots \geq \ga_n'$ such that $\ga_i=\gb+\ga_i'$, for all $1\leq i\leq n$. Hence, $\ga=\gw^{\gb}\ga'$, where $\ga'=\gw^{\ga_1'}+\cdots + \gw^{\ga_n'}$. So, $\ell(\ga)\geq \gb$ implies $\ga=\gw^{\gb+1}\ga'$, for some $\ga'>0$.  The opposite implication also holds, since $\ell(\gw^{\gb}\ga')=\gb+\ell(\ga')$.
\ep

\bcor $\ga_1+\ga_2$ is a multiple of $\gw^\gb$ iff $\ga_2$ is a multiple of $\gw^\gb$.
\ecor

\begin{lemma}[synchronization] Assume $X$, $Y$ are u.p.\ and $|X|=|Y|$. Then $X$ and $Y$ synchronize, that is, they can be represented as $X=CA^\ga$, $Y=DB^\gb$ with $|C|=|D|$, $|A|=|B|$ and $\ga=\gb$.
\end{lemma}

\bp\ By Lemma~\ref{nfperiod} we can assume $|A|=\omega^{\alpha_1}n$, $|B|=\omega^{\beta_1}m$. Consider three cases.

\medskip
\textsc{Case 1:} $\alpha_1<\beta_1$, say $\ga_1+\gy=\gb_1$ with $\gy>0$. Then  $\gw^{\gb_1}=\gw^{\ga_1}\gw^\gy =\gw^{\ga_1}n\gw^\gy$. Hence, for some $\gd\in\Lim$, $\gw^{\gb_1}=|A|\cdot \gd=|A^\gd|$.
\ignore{ Then  $\gw^{\gb_1}m=\gw^{\ga_1}\gw^\gy m=\gw^{\ga_1}n\gw^\gy m$. Hence, for some $\gd\in\Lim$, $|B|=|A|\cdot \gd=|A^\gd|$. }

We can select $D$ at least as long as $C$. Then $|DB|=|C|+\mu$, for some $\mu=\nu+\gw^{\gb_1}$  in Cantor normal form. Such a $\mu$ is a multiple of $\gw^{\gb_1}$, hence $\mu=|A^{\gy_1}|$, for some $\gy_1$. It follows that $|DB|=|CA^{\gy_1}|$. This means that $X$ and $Y$ synchronize with offset length $|DB|$ and period length $|B|$.

\medskip
\textsc{Case 2:} $\beta_1=\alpha_1$. There are two subcases.

\textsc{Case 2.1:} $\ga>\gw$. We may assume $|C|\geq |D|$. We claim that $X$ and $Y$ synchronize with offset length $|CA^\gw|$ and period length $\omega^{\alpha_1}\cdot \mathrm{lcm}(n,m)$. Since $\ga>\gw$ the word $CA^\gw$ is a proper initial segment of $X$.

Similarly to \textsc{Case 1}, we have $|CA^\gw|=|D|+\nu+\gw^{\ga_1+1}$, for some $\nu$. The ordinal $\mu=\nu+\gw^{\ga_1+1}$ is a multiple of $\gw^{\ga_1+1}$, hence $\mu=|B^\gl|$, for some $\gl\in\Lim$. It follows that $|CA^\gw|=|DB^\gl|$. Clearly, $\omega^{\alpha_1}\cdot \mathrm{lcm}(n,m)$ is a multiple of both $|A|$ and $|B|$, hence a common period length.

\medskip
\textsc{Case 2.2:} $\ga=\gw$.  We may select offsets $C$ and $D$ such that $|C|$ and $|D|$ are multiples of $\gw^{\ga_1}$. Indeed, $|CA|=\gw^{\ga_k}+\cdots+\gw^{\ga_2}+\gw^{\ga_1}$ in Cantor normal form, in particular $\ga_i\geq \ga_1$, for all $i$. Then $|CA|=\gw^{\ga_1}\cdot \ga'$, for some $\ga'>0$. Since $|CA|<|X|=\gw^{\ga_1+1}$ we have $\ga'<\gw$.

 Furthermore, we can assume $|C|\leq |D|<|CA|$. Let $A=A_1A_2$ with $|D|=|CA_1|$, then $X=CA_1(A_2A_1)^\gw= DB^\gw$. Since $|D|$ is a multiple of $\gw^{\ga_1}$ so is $|A_1|$. Since $|A|$ is a multiple of $\gw^{\ga_1}$ so is $|A_2|$. Then we have $|A_2A_1|=|A_1A_2|=|A|$ and we can select the common period length $\omega^{\alpha_1}\cdot \mathrm{lcm}(n,m)$ again.

\medskip
\textsc{Case 3:} $\beta_1<\alpha_1$. This is symmetric to \textsc{Case 1.}
\ep

\bcor \label{boo} The set $\{X\subseteq \Omega: \text{$X$ is u.p.\ in $\Omega$}\}$ is closed under the boolean operations.
\ecor

\bp\ The closure under complement is easy: the inversion of an u.p.\  $\Omega$-word is ultimately periodic.  To show the closure under union first represent the two sets $X,Y$ as synchronized $\Omega$-words. Then $X\cup Y$ corresponds to the bitwise maximum of the two words. Its period is the bitwise maximum of the corresponding periods. \ep

\subsection{Hereditarily periodic sets}

$A$ is \emph{hereditarily periodic} (h.p.)\ in $\Omega$ if, for any limit ordinal $\beta\leq\Omega$, $A$ is u.p.\ in $\beta$. Let $\HH_0(\Omega)$ denote the set of all $X\subseteq\Omega$ such that $X$ is h.p.\ in $\Omega$. We summarize basic closure properties of h.p.\ sets in the following lemma.

\bl
\benr
\item $\HH_0(\gb)$ is closed under boolean operations;
\item Any interval contained in $\gb$ is h.p.\ in $\gb$;
\item If $A\in \HH_0(\gb)$ then $(A\cap \ga) \in \HH_0(\ga)$, for all $\ga<\gb$;
\item $A\in \HH_0(\gb)$ iff $\ga+A\in \HH_0(\ga+\gb)$;
\item If $A\in\HH_0(\ga)$ and $B\in \HH_0(\gb)$ then $AB\in \HH_0(\ga+\gb)$;
\item If $A\in \HH_0(\gb)$ then $A^\ga\in \HH_0(\gb\ga)$ for all $\ga>0$.
\eenr
\el
\bp\ Claim (i) follows from Corollary~\ref{boo}. Claims (ii) and (iii) are obvious. Claim (iv) holds because the map $x\mapsto \ga + x$ preserves ultimate periodicity. Claim (v) follows from (iv), since the concatenation $AB$ is the union of sets $A$ and $|A|+B$.
Claim (vi) is proved by transfinite induction on $\ga$ using (v) for the successor step. \ep

\ignore{
Let $\gl\leq\gb\ga$ be a limit ordinal. There is a $\nu\leq\ga$ and $\gy<\gb$ such that $\gl=\gb\nu+\gy$. If $\gy>0$ then $\gy\in\Lim$ and $A^\ga\cap\gl=A^{\nu}A_1$ where $A_1=A\cap\gy$. Since $A\in\HH_0(\gb)$  $A_1$ is u.p.\ in $\gy$, hence $A^\nu A_1$ is u.p.\ in $\gl$.

If $\gy=0$ then $A^\ga\cap\gl=A^{\nu}$. If $\nu\in\Lim$ then $A^\nu$ is u.p.\ in $\gl$. If $\nu=\nu'+1$ then $A^\nu=A^{\nu'}A$ with $A$ u.p.\ in $\gb$. Hence, $A^\nu$ is u.p.\ in $\gl$. \ep }

Let $d_0$ denote the derivative operation w.r.t.\ the interval topology on $\Omega$.

\bpr \label{d0H0} $\HH_0(\Omega)$ is closed under $d_0$. \epr

\bp\ It is sufficient to show that if $X$ is u.p.\ in an $\Omega\in\Lim$ then so is $d_0(X)$.
Let $X=BA^\ga$ with $\ga\in \Lim$ and $|X|=\Omega$. We can assume $\ga$ to be a power of $\omega$. Let $\gb:=|A|$ and $\gl:=\gb\ga=|A^\ga|$. By the distributivity of $d_0$ over finite unions it is sufficient to show that $d_0(A^\ga)$ is u.p.\ in $\gl$. ($d_0(B)$ is bounded in $\Omega$, hence u.p.\ in $\Omega$.) Let $E:=\{\gb\gy:\gy<\ga\}$.

If $A=\emptyset$, then $d_0(A^\ga)=\emptyset$ and there is nothing to prove.\footnote{Pay attention that $A=\emptyset$ corresponds to a word of length $\gb$, not to the empty word.}   Assume $A\neq\emptyset$. If $\gb\in d_0(A)$ then $E\setminus\{0\}\subseteq d_0(A)$. Letting $C:=d_0(A)\cap \gb$ and $C_0:= C\cup \{0\}$ we obtain
$d_0(A^\ga)=CC_0^\ga$ which is u.p.\ in $\gl$.

If $\gb\notin d_0(A)$ then we let $E':=d_0(E)$ and observe that $d_0(A^\ga)=C^\ga\cup E'$. If $\ga=\gw$ then $E'$ is empty. Otherwise, $\ga=\gw^{1+\ga'}$ with $\ga'>0$. Then $E'=\{\gb\gw\gy:\gw\gy<\ga, \gy\neq 0\}$ is u.p.\ with period $\gb\gw$ and exponent $\gw^{\ga'}$. \ep
\bcor
$(\HH_0(\Omega),d_0)$ is a Magari algebra, hence its logic contains $\GL$.
\ecor

In the following we will show that the logic of this algebra is exactly $\GL$ provided $\Omega\geq \gw^\gw$. However, for the completeness proof for $\GLP$ this algebra of sets is insufficient. Therefore, in the next section we introduce a hierarchy of more general classes called \emph{$n$-periodic sets}.
\section{$n$-periodic sets}

\bd Define a sequence of subsets $\VV_n(\Omega)$ of $\mathcal{P}(\Omega)$ by $$\VV_0(\Omega):=\HH_0(\Omega),\quad \VV_{n+1}(\Omega):=\{B\cap \ell^{-1}(A):B\in\HH_0(\Omega), A\in\VV_n(\Omega)\}.$$
Let $\HH_n(\Omega)$ denote the closure of $\VV_n(\Omega)$ under finite unions. We call sets in $\HH_n(\Omega)$ \emph{$n$-periodic in $\Omega$}.
\ed
Intuitively, $(n+1)$-periodic sets are generated by $n$-periodic sets and all sets whose ranks are $n$-periodic. The sets $\VV_n(\Omega)$ give a more convenient base for $n$-periodic sets.

\bl For each $\Omega$,
\benr
\item $\VV_n(\Omega)$ is closed under $\cap$;
\item $\VV_n(\Omega)\subseteq \VV_{n+1}(\Omega)$;
\item $A\in \HH_n(\Omega)$ implies $\ell^{-1}(A)\in\HH_{n+1}(\Omega)$;
\item $\HH_n(\Omega)$ is closed under Boolean operations.
\eenr
\el

For a class of sets $\HH$, let $\HH|_\Omega$ denote the class  $\{A\cap\Omega:A\in\HH\}$.

\bl\ If $\Omega_1<\Omega_2$ then  $\VV_n(\Omega_1)=\VV_n(\Omega_2)|_{\Omega_1}$ and $\HH_n(\Omega_1)=\HH_n(\Omega_2)|_{\Omega_1}.$
\el
\bp\ Induction on $n$. If $n=0$ then the claim is obvious by the property of hereditarity. For $A,B\in \VV_n(\Omega_2)$ we have that $$B\cap \ell^{-1}(A)\cap\Omega_1=(B\cap\Omega_1)\cap \ell^{-1}(A\cap \Omega_1)\in \VV_{n+1}(\Omega_1).$$ Hence, $\VV_{n+1}(\Omega_2)|_{\Omega_1}\subseteq \VV_{n+1}(\Omega_1)$.
\ep
In view of this lemma, we can consider $\HH_n(\Omega)$ as the restriction of a fixed class-size algebra $\HH_n$ to $\Omega$. This allows us to write $\HH_n$ for $\HH_n(\Omega)$ whenever $\Omega$ is understood from context.

\ignore{
\bd
Let $\tau_n$ denote Icard's topology $\vartheta_{n+1}$. In other words, $\tau_0$ is the interval topology on $\Omega$, and $\tau_n$ is generated by $\tau_0$ and all sets $U^k_\gb(\Omega):=\ell^{-k}(\gb,\Omega)\cap \Omega$ where $k\leq n$ and $\gb<\Omega$. Let $d_n$ denote the derivative operation w.r.t.\ $\tau_n$ on $\Omega$.
\ed
}

Recall that $d_n$ denotes the derivative operation w.r.t.\ Icard's topology $\vartheta_n$ on $\Omega$.
Our next goal is to show that $\HH_n$ is closed under $d_0,\dots,d_n$, for $n>0$. As a base case we are going to prove that $d_n(A)\in\HH_n$, for any $A\in\HH_0$. We need a few auxiliary lemmas.

\bl For any set $A\subseteq \Omega$, $d_0(\ell^{-1}(A))=\ell^{-1}(\min A,\Omega)$ in $(\Omega,\vartheta_0)$. \el

\bp\ Recall that $\ell:(\Omega,\vartheta_0)\to (\Omega,\tau_{\leftarrow})$ is the rank function of $\vartheta_0$, hence a $d$-map. Hence, $d_0(\ell^{-1}(A))=\ell^{-1}(d_{\leftarrow}(A))=\ell^{-1}(\min A,\Omega)$. \ep

Recall that $\ell^{-1}(\min A,\Omega)$ is the set $U^1_{\beta}(\Omega)$, where $\gb=\min A$.  Corollary~\ref{Uper} below shows that sets $U^1_\gb(\Omega)$ are h.p.\ in $\Omega$.

If $A$ is a class of ordinals, we use notation $\ga\cdot A$ for $\{\ga\gb:\gb\in A\}$.

\bl \label{mult0}
$\gw^{\gb+1}\cdot [1,\gl) = U^1_\gb(\gw^{\gb+1}\gl)$.
\el

\ignore{
\bl \label{mult0}
$\gw^{\gb+1}\cdot [1,\nu]=[1,\gw^{\gb+1}\nu]\cap U^1_\gb$.
\el
}

\bp\ The inclusion $(\subseteq)$ clearly holds: If $\ga\in [1,\gl)$ then $\ell(\gw^{\gb+1}\cdot\ga)=\gb+1+\ell(\ga)\geq \gb+1>\gb$.
 For the inclusion $(\supseteq)$ write $\ga\in U^1_\gb$ in Cantor normal form
$\ga=\gw^{\ga_1}+\cdots +\gw^{\ga_n}$. We have $\ga_1\geq \cdots\geq \ga_n$. Since $\ell(\ga)=\ga_n>\gb$ there are $\ga_1'\geq \cdots \geq \ga_n'$ such that $\ga_i=\gb+1+\ga_i'$, for all $1\leq i\leq n$. Hence, $\ga=\gw^{\gb+1}\ga'$, where $\ga'=\gw^{\ga_1'}+\cdots + \gw^{\ga_n'}$. So, $\ell(\ga)>\gb$ implies $\ga=\gw^{\gb+1}\ga'$, for some $\ga'>0$. Moreover, $\ga'<\gl$, otherwise $\ga=\gw^{\gb+1}\ga'\geq\gw^{\gb+1}\gl$.
\ep

\bcor  \label{Uper} For any $\gb$, $U^1_\gb$ is hereditarily periodic in any limit ordinal.
\ecor

\bp\ Let $\gl$ be a limit ordinal. The set $U^1_\gb:=\gw^{\gb+1}\cdot [1,\gl)$ can be represented as the word $BA^\gl$ where $B$ corresponds to  the empty subset of $\gw^{\gb+1}$ and $A\subseteq \gw^{\gb+1}$ is the singleton $\{0\}$. Both $A$ and $B$ are h.p.\ in $\gw^{\gb+1}$. Hence, $U^1_\gb$ is h.p.\ in $\gw^{\gb+1}\gl$. It follows that it is also h.p.\ in $\gl$.
\ep

Let $E$ denote the set $\{\gb\ga:\ga\in\Lim\cup\{0\}\}\cap\Omega$. We note that $E\in\HH_0(\Omega)$.

\bl\ \label{dlem}
\benr
\item $d_0(E)=U^1_{\nu+1}(\Omega)$, for some $\nu$;
\item $d_n(E)= U^{n+1}_0(\Omega)\cap E$, if $n>0$;
\item For any $A\subseteq\Omega$, $d_0(E\cap \ell^{-1}(A))= U^1_{\mu}(\Omega)$, for some $\mu>0$.
\eenr
\el
Notice that the mentioned sets can also be empty.

\bp\ Since $\ga\in\Lim\cup\{0\}$ we have (letting $\ga=\gw\ga'$)
$$E = \{\gb\gw\ga':\ga'\in\On\}\cap \Omega.$$
The ordinal $\gb\gw$ has the form $\gw^{\nu+1}$, for some $\nu$.
Hence, $E=U^1_\nu(\Omega)$.

(i) $d_0(E)=d_0(U^1_\nu(\Omega))=U^1_{\nu+1}(\Omega)$.

(ii) Assume $n>0$, then $E=U^1_\nu(\Omega)$ is open in $\vartheta_n$. Therefore, $d_n(E)=d_n(\Omega)\cap E$. Since $\ell^{n+1}$ is the rank function of $\vartheta_n$ we obtain $d_n(\Omega)=\ell^{-(n+1)}(0,\Omega)=U^{n+1}_0(\Omega),$ hence the result.

(iii) We have $$U^1_{\nu}(\Omega)\cap \ell^{-1}(A)=\ell^{-1}(\nu,\Omega)\cap \ell^{-1}(A)=\ell^{-1}(A\cap (\nu,\Omega)).$$
Hence, $d_0(E)=d_0(\ell^{-1}(A\cap (\nu,\Omega)))=\ell^{-1}(\min(A\cap (\nu,\Omega)),\Omega)=U^1_\mu(\Omega)$ with $\mu= \min(A\cap (\nu,\Omega))$.
\ep

The following lemma strengthens Proposition~\ref{d0H0}. Its proof follows the same pattern, however we will refer to Proposition~\ref{d0H0} to simplify some arguments in the proof.

\bl Suppose $X$ is u.p.\ in $\Omega$. Then so is $d_0(X\cap \ell^{-1}(A))$.
\el

\bp\ Let $B$ be the period of $X$, so that $X=B_0B^\ga$, for some $\ga\in\Lim$ (which we will also assume to be a power of $\gw$). Let $\ga_0=|B_0|$ and $\gb:=|B|=\gw^\nu(n+1)$.

The subset $X$ is naturally split in two subsets of $\Omega$  corresponding to the subwords $B_0$ and $B^\ga$. Hence, by the additivity of $d_0$ it is sufficient to show that $d_0(B_0\cap\ell^{-1}(A))$ u.p.\ in $\Omega$ and $d_0(B^{\ga}\cap\ell^{-1}(A))$ is u.p.\ in $\gl:=\gb\ga$. (The bound $\Omega$ changes to $\gb\ga$, since by considering $B^\ga$ we shift the set leftwards by $\ga_0$. This operation does not change the rank\footnote{That is, the value of the $\ell$-function.} of any point except for the first bit in $B^\ga$. Changing a set at one point has no effect on the derivative.)

The set $B_0\cap\ell^{-1}(A)$ is bounded in $\Omega$, hence so is $d_0(B_0\cap\ell^{-1}(A))$, and therefore it is u.p.\ in $\Omega$.

To analyse $d_0(B^\ga\cap\ell^{-1}(A))$ we consider the set  $E:=\{\gb\gy : \gy<\ga\}.$ These are the coordinates of the first bits of occurrences of period $B$.
We split $B^\ga$ into sets $B^\ga\cap E$ and $B^\ga\setminus E$.
By periodicity, if $0\in B$ then $E\cap B^\ga=E$, otherwise $E\cap B^\ga=\emptyset$.

Using the additivity of $d_0$ we now deal separately with the sets $d_0((B^\ga\setminus E)\cap \ell^{-1}(A))$ and $d_0((B^\ga\cap E)\cap \ell^{-1}(A))$ on $\gl$.

For the former, we observe that the ranks of points in $B^\ga\setminus E$ are the same as of the corresponding points in $B\setminus\{0\}\subseteq \gb$. Hence, $d_0((B^\ga\setminus E)\cap \ell^{-1}(A))=d_0(C^\ga),$ where $C:=(B\setminus\{0\})\cap \ell^{-1}(A)$ with $|C|=\gb$.
Therefore, this set is the derivative of an u.p.\ set in $\gl$, hence u.p.\ itself by Proposition~\ref{d0H0}.

For the latter, we consider two cases.  If $0\notin B$ then $d_0((B^\ga\cap E)\cap \ell^{-1}(A))$ is empty. If $0\in B$ then $d_0((B^\ga\cap E)\cap \ell^{-1}(A))=d_0(E\cap \ell^{-1}(A))$.
Here we consider two subcases according to the value of the exponent $\ga$ (which is a power of $\gw$). If $\ga=\gw$ then $d_0(E\cap \ell^{-1}(A))\cap \gl\subseteq d_0(E)\cap\gl=\emptyset$. If $\ga=\gw^{\nu}$ with $\nu>1$, then by Lemma~\ref{dlem} (iii) we obtain $d_0(E\cap \ell^{-1}(A))=U^1_\mu(\gl)$, for some $\mu$. Hence, it is also u.p.\ in $\gl$.~\ep

\ignore{
Case 1: $(B_1\setminus \{0\})\cap \ell^{-1}(A)=\emptyset$. If $0\notin B_1$, then $B_1^\ga \cap \ell^{-1}(A) = \emptyset$.

\ignore{ and $B\cap \ell^{-1}(A)=B_0\cap \ell^{-1}(A)$ is bounded in $\Omega$. Hence, so is $d_0(B\cap\ell^{-1}(A))$, and therefore it is ultimately periodic in $\Omega$.}

If $0\in B_1$ then $E\subseteq B_1^\ga$ and $B_1^\ga\cap \ell^{-1}(A)=E\cap  \ell^{-1}(A)$. By Corollary~\ref{E-cor},  $d_0(E\cap  \ell^{-1}(A))$ is of the form $U^1_\gb$, for some $\gb$.

Case 2: $(B_1\setminus\{0\})\cap \ell^{-1}(A)\neq \emptyset$. Let $C:=d_0(B_1\cap \ell^{-1}(A))\cap \ga_1$ considered as a word of length $|C|=\ga_1$. Clearly, $C^\ga\subseteq d_0(B_1^\ga\cap \ell^{-1}(A))$. We need to determine which points of $E$ belong to $d_0(B_1^\ga\cap \ell^{-1}(A))$.

Subcase 2.1: $B_1\cap \ell^{-1}(A)$ is unbounded in $\ga_1$. Then $E\setminus\{0\}\subseteq d_0(B_1^\ga\cap \ell^{-1}(A))$. So we let $C_1:=\{0\}\cup C$ and obtain $d_0(B_1^\ga\cap \ell^{-1}(A))=CC_1^\ga$.

Subcase 2.2: $B_1\cap \ell^{-1}(A)$ is bounded in $\ga_1$. Then $d_0(B_1^\ga\cap \ell^{-1}(A))=C^\ga\cup E'$ where $E'=d_0(E)$ is the set of all limit points in the natural enumeration of $E$. If $\ga=\gw$ then $E'$ is empty. Otherwise,  $E'=\{\ga_1\gw\gy:\gw\gy<\ga\}$ is u.p.\ with period $\ga_1\gw$.
}

\bcor \label{d0H} If $B\in\HH_0(\Omega)$ then $d_0(B\cap \ell^{-1}(A))\in \HH_0(\Omega)$.
\ecor

\bcor \label{H0n} If $V\in\HH_n$ then $d_0(V)\in\HH_0$.
\ecor

\bp\ It is sufficient to prove the lemma for $V\in \VV_n$. If $n=0$ the statement is in Lemma \ref{d0H0}. If $n>0$, then $V=B\cap \ell^{-1}(A)$, for some $B\in\HH_0$ and $A\in\VV_{n-1}$. By Lemma \ref{d0H} $d_0(V)\in\HH_0$.  \ep

\bl \label{ldist} If $d$ is a derivative operator of a topology $\tau$ extending $\vartheta_1$, then, for all $B\in\HH_0$,
$$d(B\cap \ell^{-1}(A))=d(B)\cap d(\ell^{-1}(A)).$$
\el

\bp\ Let $\gl\leq\Omega$, $\gl\in d(B)\cap d(\ell^{-1}(A)).$ Then $B$ is u.p.\ in $\gl$ that is
$B=B_0B_1^\ga$ with $\ga\in\Lim$. As before, without loss of generality we may assume $B_0$ to be empty.

Since $\gl\in d(B)\subseteq d_1(B)$ we have $\ell^2(\gl)>0$ and we can assume that $|B_1|=\gw^\nu$ with $\nu>0$.

As before, we let $E:=\{\gw^{\nu}\gy:\gy<\ga\}$. If $0\in B_1$ then $E\subseteq B_1^\ga$, otherwise $B_1^\ga\cap E=\emptyset$. %so $B_1^\ga=(B_1\setminus\{0\})^\ga\cup E$. Otherwise, $B_1^\ga=(B_1\setminus\{0\})^\ga$.

For any point $x<\gw^\nu$ we have $\ell(x)<\nu$, whereas $\ell(\gl)=\nu+\ga>\nu$. By the periodicity and preservation of ranks it follows that $\ell(x)<\nu$ for all $x\in \gl\setminus E$.

The set $U^1_\nu=\ell^{-1}(\nu,\gl)\cap\gl$ is a $\tau$-open punctured neighbourhood of $\gl$. By the above, $U^1_\nu\subseteq E$. To show that $\gl\in d(B_1^\ga\cap \ell^{-1}A)$ consider any $\tau$-open punctured neighbourhood $U$ of $\gl$. We may assume that $U\subseteq U^1_\nu$, otherwise consider their intersection.

Since $\gl\in d(\ell^{-1}A)$ there is a point $x\in U\cap \ell^{-1}A$. Then $x\in E$ and $\ell(x)\in A$.
On the other hand, $\gl\in d(B_1^\ga)$, hence there is a point $y\in U\subseteq E$ such that $y\in B_1^\ga$. Therefore, $0\in B_1$ and $E\subseteq B_1^\ga$, so $x\in B_1^\ga\cap \ell^{-1}A$ and $x\in U$, as required.
\ep

\bl \label{dnH0} If $X\in \HH_0(\Omega)$ and $n>0$, then $d_n(X)\in \HH_n(\Omega)$.
\el
\bp\ The proof generalizes that of Lemma \ref{d0H0}. We argue by transfinite induction on $\Omega$. So, we assume the claim to be true for all limit $\gl<\Omega$.
 %It is sufficient to show that if $X$ is u.p.\ in $\Omega$ then $d_n(X)\in \HH_n(\Omega)$.

Let $X=BA^\ga$ with $\ga\in \Lim$ and $|X|=\Omega$. Let $\gb:=|A|$ and $\gl:=\gb\ga=|A^\ga|$. By the distributivity of $d_n$ over finite unions it is sufficient to show that $d_n(B)\in\HH_n(\Omega)$ and $d_n(A^\ga)\in\HH_n(\gl)$. The former is equivalent to $d_n(B)\in\HH_n(\mu)$ for $\mu:=|B|$. Since $\mu<\Omega$ we have the claim by the transfinite induction hypothesis. To deal with the latter we consider two sets.

Let $I_\gy$ denote the $\vartheta_0$-clopen interval $[\gb\gy+1,\gb(\gy+1)]$. Then $\gl$ is the disjoint union of sets $I:=\bigcup_{\gy<\ga}I_\gy$ and $E:=\{\gb\gy:\gy<\ga, \gy\in\Lim\cup\{0\}\}$.
%Given any point $x\in\gl$ we determine if $x\in d_n(A^\ga)$ depending on whether $x\in I$ or $x\in E$.
We have $d_n(A^\ga)\cap\gl=(d_n(A^\ga)\cap I)\cup (d_n(A^\ga)\cap E)$, and we consider the two parts $I$ and $E$ separately.

\bs
\textsc{Case $I$.} Since $I_\gy$ is clopen, $d_n(A^\ga)\cap I_\gy=d_n(A^\ga\cap I_\gy)$. By periodicity, the latter is isomorphic to $C:=d_n(A)\cap I_0$ by the map $f_\gy:I_0\to I_\gy$, $x\mapsto \gb\gy+x$. So, $$d_n(A^\ga)\cap I=\bigcup_{\gy<\ga}f_\gy(C).$$

Let $C_0:= (C\cap \gb)\cup\{0\}$ if $\gb\in C$, and $C_0:=C\cap \gb$, otherwise. Then $d_n(A^\ga)\cap I=C_0^\ga\cap I$. As both $E,I\in\HH_0(\gl)$ it is now sufficient to prove that $C_0^\ga\in\HH_n(\gl)$.
If $n=0$ and $C_0\in\HH_0(\gl)$, then $C_0^\alpha$ is u.p.\ in any ordinal $\mu\leq\gl$. If $n>0$ then $C_0$ has the form $\bigcup_{i=1}^k (B_i\cap \ell^{-1}A_i)$ with $B_i\in \HH_0(\gb)$ and $A_i\in \VV_{n-1}(\gb)$. Notice that $f_\gy$ preserves the ranks of points, hence $$f_\gy(B_i\cap \ell^{-1}(A_i))=f_\gy(B_i)\cap \ell^{-1}(A_i).$$ Then $$\bigcup_{\gy<\ga}f_\gy(C)=\bigcup_{\gy<\ga}\bigcup_{i=1}^k (f_\gy(B_i)\cap \ell^{-1}(A_i))=\bigcup_{i=1}^k(\bigcup_{\gy<\ga}f_\gy(B_i)\cap \ell^{-1}(A_i)). $$
Letting $D_i:= \bigcup_{\gy<\ga}f_\gy(B_i)$ thus we obtain $$d_n(A^\ga)\cap I= \bigcup_{i=1}^k(D_i\cap \ell^{-1}(A_i)).$$

The sets $D_i\in \HH_0(\gl)$ by the same argument as for the case $n=0$ (with $D_i$ playing the role of $C$). Hence, $D_i\cap \ell^{-1}(A_i)\in\VV_{n}(\gl)$ and then $d_n(A^\ga)\cap I\in \HH_n(\gl)$.

\bs
\textsc{Case $E$.}
Since $n>0$ the set $E$ is clopen in $\vartheta_n$. Therefore,  $d_n(A^\ga)\cap E=d_n(A^\ga\cap E)$. However, $A^\ga\cap E$ is empty, if $0\notin A$, and $A^\ga\cap E = E$, otherwise. In the first case $d_n(A^\ga\cap E)$ is empty. In the second case, Lemma~\ref{dlem} (ii) is applicable and we obtain $d_n(A^\ga\cap E)=d_n(E)=U_0^{n+1}(\gl)\cap E$. We have $U_0^{n+1}(\gl)=\ell^{-n}(\gl\cap\Lim)$. Since $\gl\cap\Lim\in\HH_0(\gl)$ we conclude that $U_0^{n+1}$ and $d_n(E)$ are in $\HH_n(\gl)$, as required. \ep

\bl If $V\in \HH_n$, $n>0$, then $d_i(V)\in \HH_i$ for all $i\leq n$. \el

\bp\ It is sufficient to prove it for $V\in\VV_n$. We use induction on $i\leq n$. For $i=0$ the statement follows from Lemma \ref{H0n}. For $i>0$ we have $d_i(V)=d_i(B\cap \ell^{-1}(A))$ with $B\in\HH_0$ and $A\in\VV_{n-1}$. Then by Lemma \ref{ldist} $d_i(V)=d_i(B)\cap d_i(\ell^{-1}(A))=d_i(B)\cap \ell^{-1}(d_{i-1}A)$. By Lemma~\ref{dnH0} $d_i(B)\in \HH_i$ and by the induction hypothesis $d_{i-1}A\in \HH_{i-1}$, hence $\ell^{-1}(d_{i-1}A)\in \HH_{i}$.
\ep

\section{Logic $J_n$ and $J_n$-morphisms}

This section presents the technical notion of $J_n$-morphism from~\cite{BekGab13}. This notion allows us to map general topological frames to Kripke models of logic $J_n$ in such a way that the validity of formulas is preserved. It is essentially the same as the notion of $J_n$-morphism from~\cite{BekGab13} with the exception that we now also require backwards preservation of admissible sets. This requirement is relatively minor, so our presentation is almost literally taken from~\cite{BekGab13} and all the results in this section are from that paper.

Let $\mathcal{L}_n$ denote the modal language with modalities
$[0],[1],\dots,[n]$. Recall that we denote by $\GLP_n$ the restriction of $\GLP$ to $\mathcal{L}_n$.
While $\GLP$ and $\GLP_n$ for $n>0$ are Kripke incomplete, the following subsystem $\J$ serves as a best approximation to $\GLP$ which is. An axiomatization of $\J$ is obtained by replacing Axiom (\ref{mon}) by two theorems of $\GLP$:

\benr
\item $ [m]\phi
\imp [n][m]\phi$, for $m\leq n$;
\item $[m]\phi\to [m][n]\phi$, for $m\leq n$.
\eenr

We denote by $\J_n$ the fragment of $\J$ in $\cL_n$.

\bd
The class of \emph{$J_n$-trees}, which are finite Kripke frames of the form $(W,R_0,\dots,R_n)$, is defined inductively as follows.
\ben
\item $J_0$-trees are finite irreflexive transitive tree-like Kripke frames $(W,R_0)$;
\item If $\cW=(W,R_0,\dots,R_n)$ is a $J_n$-tree, then $\cW^+:=(W,\emptyset,R_0,\dots,R_n)$ is a $J_{n+1}$-tree;
\item If $\cW_1, \dots, \cW_k$ are $J_n$-trees and $\cV$ is a $J_{n-1}$-tree, then the structure $\cW=(W,R_0,\dots,R_n)$ is a $J_n$-tree, where $\cW$ is obtained from the disjoint sum of structures $\cV^+$ and $\cW_1,\dots,\cW_k$ by extending $R_0$ with all pairs $(v,w)$ for $v\in V$ and $w\in \bigcup_{i=1}^k  W_i$.
\een
We call \emph{$J$-trees} the structures that are $J_n$-trees for some $n\in\gw$. We can view them as Kripke frames $(W,(R_i)_{i\in\gw})$ where all relations $R_i$ for $i>n$ are empty.
\ed

It is easy to check that the relations of a $J$-tree are irreflexive partial orders and satisfy the following properties:
\bi \item $R_kR_m\subseteq R_{\min(k,m)}$ (polytransitivity); \item $R_m^{-1}R_k\subseteq R_k$, for $k<m$.\ei

These conditions ensure that $\J_n$ is sound w.r.t.\ $J_n$-trees. The following proposition is from~\cite{Bek10}.
\bpr \label{jlem}
\benr
\item $\J_n$ is complete w.r.t.\  the class of $J_n$-trees;
\item $\J$ is complete w.r.t.\ the class of $J$-trees.
\eenr
\epr

Let $\cT=(T,R_0,\dots, R_n)$ be a $J_n$-tree and let $R^*_k$ denote $\bigcup_{i=k}^n R_i$. By the polytransitivity, $R^*_k$ is transitive.
An $R^*_k$-minimal node $w\in T$ is called a \emph{hereditary  $k$-root}. Since $\cT$ is a $J_n$-tree, for each $w\in T$ and each $k\leq n$ there exists a hereditary $k$-root $v\in T$ such that $v=w$ or
$vR_k w$. We call the \emph{root} of $\cT$ its hereditary $0$-root.

As any Kripke model, we can view $\cT$ as a poly-topological space $\cT=(T,\sigma_0,\dots,\sigma_n)$ by considering all $R_i$-upsets to be $\sigma_i$-open.

\bd
Given a poly-topological space $\cX=(X,\tau_0,\dots,\tau_n)$ and a map $f:\cX\to \cT$ we say that $f$ is a \emph{$J_n$-morphism} if:

\begin{itemize}
%\item[($j_0$)] $f^{-1}(A)\in \HH$ for all $A\subseteq T$;
\item[($j_1$)] $f:(X,\tau_n)\to (T,\sigma_n)$ is a $d$-map;
\item[($j_2$)] $f:(X,\tau_k)\to (T,\sigma_k)$ is an open map for all $k\leq n$;
\item[($j_3$)] For each $k<n$ and each hereditary $(k+1)$-root $w\in T$, the sets $f^{-1}(R_k^*(w))$ and $f^{-1}(R_k^*(w)\cup\{w\})$ are open in $\tau_k$;
\item[($j_4$)] For each $k<n$ and each hereditary $(k+1)$-root $w\in T$, the set $f^{-1}(w)$ is a $\tau_k$-discrete subspace of $\cX$.
\end{itemize}

Here $R_k^*(w)$ denotes the set $\{u\in T: w R_k^* u\}$. Notice that the sets $R_k^*(w)$ and $R_k^*(w)\cup \{w\}$ are $\gs_i$-open for all $i\geq k$, therefore condition ($j_3$) is a weakening of the continuity of $f:(X,\tau_k)\to (T,\gs_k)$.

If $f:\cX\to \cT$ is a $J_n$-morphism and $(\cX,\cA)$ is an ambience on $\cX$, then $f$ is called \emph{algebra preserving} if $f^{-1}(A)\in \cA$ for all $A\subseteq T$. (We consider $\cT$ as an ambience where all subsets are admissible.) Since $T$ is finite, this requirement is equivalent to $f^{-1}(x)\in \cA$, for all $x\in T$.
\ed

Let $\phi$ be a formula in $\cL_n$. Define:  \[M(\phi) :=
\bigwedge_{i<s}\bigwedge_{k=m_i+1}^n
([m_i]\phi_i\rightarrow[k]\phi_i),\] where $[m_i]\phi_i$, $i<s$, are
all subformulas of $\phi$ of the form $[m]\psi$ and
$n:=\max_{i<s}m_i$. Also, let
$M^+(\phi):=M(\phi)\land\bigwedge_{m\leq n}[m]M(\phi)$.
These formulas provide a well-known translation from $\GLP_n$ to $\J_n$. Clearly, $\GLP_n\vdash M^+(\phi)$.

The following theorem allows one to reduce the completeness proof for $\GLP$ to a covering construction. It is obtained by the same proof as in~\cite{BekGab13} for the situation of ambiences. For the reader's convenience, we reproduce the proof of this theorem in the Appendix.

\bt \label{t:main-1}
Let $\cX$ be an ambience, $\cT$ a ${J}_n$-tree, $f:\cX\tto \cT$ an algebra preserving surjective
${J}_n$-morphism. Then, for all $\mathcal{L}_n$-formulas  $\varphi$,
$\cX\models\varphi$ implies $\cT\models M^+(\varphi)\to\varphi$.
\et

\section{Covering lemma}

Recall that $\ge_0$ is the supremum of
the countable ordinals $\gw_k$ recursively defined by $\gw_0=1$ and $\gw_{k+1}=\gw^{\gw_k}$. The crucial point of the completeness proof is the following lemma.

\bl[covering] \label{lkey}
For each ${J}_n$-tree $\cT$ there exist an ordinal $\gl<\gw_{n+2}$ and an onto algebra preserving ${J}_n$-morphism
$f:\cX\tto \cT$ where
$\cX=([1,\gl],\vartheta_0,\dots,\vartheta_n,\HH_n)$.
\el

Of course, the ambience $\cX$ is isomorphic to $(\Omega,\vartheta_0,\dots,\vartheta_n,\HH_n)$, where $\Omega=\gl+1$ if $\gl$ is infinite, otherwise $\Omega=\gl$, but it will be more convenient in the inductive construction below to deal with the intervals of the form $[1,\gl]$.
For the proof we need a few preparatory lemmas. The first one covers the case $n=0$.

\bl \label{base} For each finite tree $\cT=(T, R_0)$ with the root $r$ there is an ordinal $\gl<\gw^\gw$ and an onto algebra preserving $d$-map $f:\cX\tto \cT$ where $\cX=([1,\gl],\vartheta_0,\HH_0)$ and  $f^{-1}(r)=\{\lambda\}$.
\el

\bp\ The construction of such a $d$-map is well-known, see e.g.\ \cite{BM10}. The essentially new claim is that the $d$-map is algebra preserving. We construct $f$ by induction on the height of $\cT$. If $ht(\cT)=0$, then we put $\gl=1$. If $ht(\cT)>0$ let $\cT_1,\dots,\cT_k$ denote the immediate subtrees of $\cT$, their roots are denoted $r_1,\dots, r_k$. By the induction hypothesis there are algebra preserving $d$-maps $g_i:([1,\gl_i],\vartheta_0,\HH_0) \tto \cT_i$ such that $g_i^{-1}(r_i)=\gl_i$, $i=1,\dots, k$. We define $\mu:=\gl_1+\cdots +\gl_k$ and a $d$-map $g:[1,\mu]\tto \cS$ where $\cS$ is the disjoint sum of the trees $\cT_i$. Obviously, $([1,\mu],\vartheta_0,\HH_0)$ is isomorphic to the sum of ambiences $([1,\gl_i],\vartheta_0,\HH_0)$, so $g$ is defined as the sum of maps $f_i$ and is an algebra preserving $d$-map. Then we put $\gl:=\mu\gw$ and define $f(\gl):= r$ and $f(\ga):= g(\gb)$ if $\ga=\mu n +\gb$, where $n<\gw$ and $1\leq \gb\leq \mu$.

Obviously, $f:[1,\gl]\tto T$. The interval $[1,\gl)$ is isomorphic to the sum of $\gw$-many copies of $[1,\mu]$, and the restriction of $f$ to each copy is the same as $g$. Hence, the restriction of $f$ to $[1,\gl)$ is a $d$-map onto $\cS$.

To verify that $f$ is an open map consider an open set $U\subseteq [1,\gl]$. If $\gl\notin U$ then $f(U)$ is open by the above. Otherwise, $U$ contains an interval $[\mu n,\gl]$, for some $n$. In this case $f(U)=T$ and is also open. Similarly, it is easy to see that $f$ is continuous and pointwise discrete.

We show that $f$ is algebra preserving, that is, $f^{-1}(w)\in\HH_0$, for each $w\in T$. If $w=r$ the claim trivially holds. If $w\in T$ and $w\neq r$, then $f^{-1}(w)$ can be represented as a transfinite word $A^\gw$ where $A=g^{-1}(w)\subseteq [1,\mu]$. By the induction hypothesis $A\in\HH_0$. It follows that $f^{-1}(w)$ is hereditarily periodic in each $\ga<\gl$. On the other hand, $A^\gw$ is u.p.\ in $\gl$, hence $A^\gw$ is hereditarily periodic in $\gl$. \ep

We remark that the $d$-map constructed in Lemma~\ref{base} is also algebra preserving for each of the larger algebras $\HH_n$, for $n>0$, and $\HH$.

\bl \label{lifting}
For all $\gl$, $n$, the function $\ell$ is an algebra preserving onto $d$-map
$$\ell:([1,\gw^\gl],\vartheta_{n+1},\HH_{n+1})\to ([0,\gl],\vartheta_n,\HH_n).$$
\el
\bp\ We already know that $\ell$ is an onto $d$-map. We show that it is algebra preserving. Every set in $\HH_n$ is a finite union of sets in $\cV_n$. Since $\ell^{-1}$ preserves finite unions, it is sufficient to show that $\ell^{-1}(A)\in \cV_{n+1}$ if $A\in\cV_n$.  But $\ell^{-1}(A)\in \cV_{n+1}$ by definition.
\ep
An immediate corollary of this lemma is that $\ell$ is also an algebra preserving $d$-map  $$([1,\gw^\gl],\vartheta_{i+1},\HH_n)\to ([0,\gl],\vartheta_i,\HH_{n-1}),$$ for each $i<n$.

\proof{ of covering lemma.} For each $J_n$-tree $\cT=(T,R_0,\dots,R_n)$ with a root $r$ we are going to build a periodic frame
$\cX=([1,\gl],\vartheta_0,\dots,\vartheta_n,\HH_n)$ and an algebra preserving  onto $J_n$-morphism $f:\cX\tto \cT$. The construction goes by induction on $n$ with a subordinate induction on the $R_0$-height of $\cT$, which is denoted $ht_0(\cT)$. For the induction to work it will be convenient for us to ensure two additional properties:
$f^{-1}(r)=\{\gl\}$ and, for each $w\in T$, $f^{-1}(w)$ is a $\vartheta_0$-discrete subset of $X$. Such $J_n$-morphisms will be called \emph{suitable}.

The case $n=0$ is Lemma~\ref{base}. Assume $n>0$ and the claim
holds for all $m<n$.
If $ht_0(\cT)=0$ then $R_0=\emptyset$. Consider the $J_{n-1}$-tree $\cS:=(T,R_1,\dots,R_n)$. By the induction hypothesis, for some $\mu$, there is a suitable $J_{n-1}$-morphism  $$g:([0,\mu],\vartheta_0,\dots, \vartheta_{n-1},\HH_{n-1})\tto \cS.$$
(Here we also use the obvious isomorphism between the general frames $[0,\mu]$ and $[1,1+\mu]$.)
By Lemma~\ref{lifting}, $\ell$ is an algebra preserving onto $d$-map, for each $i<n$,
$$\ell:([1,\gw^\mu],\vartheta_{i+1},\HH_n)\tto ([0,\mu],\vartheta_i,\HH_{n-1}).$$
For this map $\ell^{-1}(\mu)=\{\gw^\mu\}$.
Hence, it is easy to verify that
$$g\circ \ell: ([1,\gw^\mu],\vartheta_0,\dots, \vartheta_{n},\HH_{n})\tto (T,\emptyset,R_1,\dots, R_n)$$ is a suitable $J_n$-morphism.

\bigskip
Assume $ht_0(\cT)>0$. Let $\cT_1,\dots, \cT_m$ denote the immediate subtrees of $\cT$ with the roots $r_1,\dots, r_m$.  By the induction hypothesis there are suitable $J_n$-morphisms $g_i:([1,\gl_i],\vartheta_0,\dots, \vartheta_{n},\HH_{n})\tto \cT_i$ such that $g_i^{-1}(r_i)=\gl_i$, $i=1,\dots, m$. We define $\nu:=\gl_1+\cdots +\gl_m$ and a map $g:[1,\nu]\tto \cS$ where $\cS$ is the disjoint sum of the trees $\cT_i$. Obviously, $([1,\nu],\vartheta_0,\dots,\vartheta_n,\HH_n)$ is isomorphic to the sum of ambiences  $([1,\gl_i],\vartheta_0,\dots,\vartheta_n,\HH_n)$, so we define $g$ as the sum of maps $f_i$. We have that $g$ is algebra preserving and satisfies the conditions ($j_1$)--($j_4$).

Further, as in base case, consider the subset $Q:=R^*_1(r)\cup \{r\}$ and the $J_{n-1}$-tree $(Q,R_1,\dots,R_n)$. Notice that the restriction of $R_0$ to $Q$ is empty. By the induction hypothesis, for some $\mu$, there is a suitable $J_{n-1}$-morphism  $$h:([1,\mu],\vartheta_0,\dots, \vartheta_{n-1},\HH_{n-1})\tto (Q,R_1,\dots,R_n).$$
By Lemma~\ref{lifting}, $\ell$ is an algebra preserving onto $d$-map, for each $i<n$,
$$\ell:([1,\gw^\mu],\vartheta_{i+1},\HH_n)\tto ([0,\mu],\vartheta_i,\HH_{n-1}).$$
%Moreover, for this map $\ell^{-1}(\mu)=\{\gw^\mu\}$.
Since $\ell^{-1}(0)$ is the set of successor ordinals in $[1,\gw^\mu]$, $\ell$ maps $Y:=[1,\gw^\mu]\cap \Lim$ onto $[1,\mu]$.
For this map $\ell^{-1}(\mu)=\{\gw^\mu\}$.
Hence, it is easy to verify that
$$h\circ \ell: (Y,\vartheta_0,\dots, \vartheta_{n},\HH_{n})\tto (T,\emptyset,R_1,\dots, R_n)$$ is a suitable $J_n$-morphism.
Then we put $\gl:=\nu\cdot\gw^\mu$ and define for $1\leq \ga\leq \gl$
$$ f(\ga) :=
\begin{cases} %r, & \text{if $\ga=\gl$}, \\
h(\ell(\gy)), & \text{if $\ga=\nu\gy$, $\gy\in\Lim$, $\gy\leq \gw^\mu$,} \\
g(\gb), & \text{if $\ga=\nu \gamma +\gb$, $\gy<\gw^\mu$, $1\leq \gb\leq \nu$.}
\end{cases}
$$

We claim that $f$ is a suitable $J_n$-morphism $([1,\gl],\vartheta_0,\dots,\vartheta_n,\HH_n) \tto \cT$.

Notice that $[1,\gl]$ is the disjoint sum of two sets: $L:=\{\nu\ga:\ga\in Y\}=\nu\cdot Y$ and $M:=\{\nu\ga+\gb:\gb\in [1,\nu],\ga<\gw^\mu\}$. Clearly, $L=f^{-1}(Q)$ and $M=f^{-1}(S)$. We also have $f^{-1}(r)=\nu\cdot \ell^{-1}(h^{-1}(r))=\{\nu\cdot\gw^\mu\},$ as expected.

The subspace $(L,\vartheta_i)$ is homeomorphic to $(Y,\vartheta_i)$ by the map $\ga\mapsto \nu\ga$. The subspace $(M,\vartheta_i)$ is the disjoint sum of clopens $I_\ga:=[\nu\ga+1,\nu\ga+\nu]$, each of which is homeomorphic to $([1,\nu],\vartheta_i)$. In particular, $M$ is $\vartheta_0$-open in $[1,\gl]$.

By Lemma~\ref{mult0} $L$ is $\vartheta_1$-open, hence $\vartheta_i$-open in $[1,\gl]$, for all $i>0$. (We have $[1,\gw^{\mu}]\cap \Lim=\{\nu\gw\ga:0<\ga\leq \gw^\mu\}$. However,  $\nu\gw=\gw^{\gb+1}$, for some $\gb$, so Lemma~\ref{mult0} applies.)

\medskip
\textsc{Claim 1:} $f$ is open for each $\vartheta_i$, $i\leq n$.

Consider a $\vartheta_i$-open set $A\subseteq [1,\gl]$. Then $f(A)=f(A\cap M)\cup f(A\cap L)$. The set $A\cap I_\ga$ is $\vartheta_i$-open, therefore so is $B:=\{\gb\in[1,\nu]: \nu\ga+\gb\in A\}$. Since $g:[1,\nu]\to S$ is a $J_n$-morphism,  $f(A\cap I_\ga)=g(B)$ is $\gs_i$-open. This holds for all $\ga<\gw^\mu$, hence $f(A\cap M)=\bigcup_{\ga<\gw^\mu}A\cap I_\ga$ is $\gs_i$-open.

If $A$ is $\vartheta_i$-open, then $A\cap L$ is $\vartheta_i$-open in the relative topology of $L$, hence $B:=\{\gy: \nu \gy\in A\}$ is open in  $([1,\gw^\mu]\cap\Lim,\vartheta_i)$. Then $\ell(B)$ is $\vartheta_{i-1}$-open and $f(A\cap L)=h(\ell(B))$ is $\gs_i$-open  in $(Q,\gs_1,\dots,\gs_n)$, as $h$ is a $J_n$-morphism from $([0,\mu],\vartheta_0,\dots,\vartheta_{n-1})$ to  $(Q,R_1,\dots,R_n)$.
We have shown that $f(A\cap M)$ and $f(A\cap L)$ are both $\gs_i$-open, hence so is $f(A)$.

 \medskip
\textsc{Claim 2:} The map $f:([1,\gl],\vartheta_n)\to (T,\gs_n)$ is continuous.% and pointwise discrete, hence satisfies $(j_1)$.

If $A\subseteq T$ is $\gs_n$-open, then so are $A\cap Q$ and $A\cap S$ ($S$ is $\gs_i$-open for all $i\leq n$, and $Q$ is $\gs_i$-open for $i>0$.) Since $h\circ\ell:Y\tto Q$ is a $J_n$-morphism, the set $B:= \ell^{-1}(h^{-1}(A\cap Q))$ is $\vartheta_n$-open in $Y$. Then $\nu \cdot B$ is $\vartheta_n$-open in $L$. Since $L$ is a $\vartheta_1$-open subset of $[1,\gl]$ and $n\geq 1$, we conclude that $f^{-1}(A\cap Q)=\nu\cdot B$ is also $\vartheta_n$-open.

Similarly,
the set $C:=g^{-1}(A\cap S)\subseteq [1,\nu]$ is $\vartheta_i$-open, since $g$ is continuous. Hence, so is the set $\nu\ga + C\subseteq I_\ga$, for each $\ga$, and therefore
$f^{-1}(A\cap S)=\bigcup_{\ga<\gw^\mu}(\nu\ga + C)$ is $\vartheta_n$-open. Thus, we have proved that $f^{-1}(A)=f^{-1}(A\cap S)\cup f^{-1}(A\cap Q)$ is $\vartheta_n$-open.

\medskip
\textsc{Claim 3:} For $k<n$, if $w$ is a $(k+1)$-root of $\cT$, then the sets $f^{-1}(R_k^*(w))$ and $f^{-1}(R_k^*(w)\cup \{w\})$ are $\vartheta_k$-open.

If $w\in S$, then $R_k^*(w)\subseteq S$, and since $g:[1,\nu]\tto S$ is a $J_n$-morphism, $g^{-1}(R_k^*(w))$ and $g^{-1}(R_k^*(w)\cup \{w\})$ are $\vartheta_k$-open in $[1,\nu]$. It follows that $f^{-1}(R_k^*(w))$ and $f^{-1}(R_k^*(w)\cup \{w\})$ are $\vartheta_k$-open in $M$ and in $[1,\gl]$, as $M$ is $\vartheta_0$-open in $[1,\gl]$.

If $w\in Q$ we consider two subcases. If $k=0$ then $w=r$ (the only 1-root of $T$ in $Q$). Then $R_0^*(w)=T\setminus\{r\}$ and hence $f^{-1}(R_0^*(w))=[1,\gl)$, which is $\vartheta_0$-open. The set $f^{-1}(R_0^*(w)\cup \{w\})=f^{-1}(T)=[1,\gl]$ is also $\vartheta_0$-open.

If $k>0$, then $R^*_k(w)\subseteq Q$. Since $h\circ \ell:[1,\gw^\mu]\tto Q$ is a $J_{n}$-morphism, we obtain that $A:=(h\circ\ell)^{-1}(R^*_k(w))$ and $B:=(h\circ\ell)^{-1}(R^*_k(w)\cup \{w\})$ are $\vartheta_k$-open in $Y$. Since $\ga\mapsto \nu\cdot\ga$ is a $\vartheta_k$-open map, for $k>0$, we conclude that $f^{-1}(R^*_k(w))=\nu\cdot A$ and $f^{-1}(R^*_k(w)\cup \{w\})=\nu\cdot B$ are $\vartheta_k$-open.

\medskip
\textsc{Claim 4}: For each $w\in T$ the set $f^{-1}(w)$ is $\vartheta_0$-discrete.

Let $w\in T$. If $w\in S$ then $f^{-1}(w)\subseteq M$ is $\vartheta_0$-discrete in the relative topology of $M$. This yields that  $f^{-1}(w)$ is $\vartheta_0$-discrete in $[1,\gl]$.
If $w\in Q$ then $h^{-1}(w)$ is $\vartheta_0$-discrete in $[1,\mu]$. By Lemma~\ref{rank-dmap} (i) $\ell$ is $\vartheta_0$-pointwise discrete, therefore   $B=\ell^{-1}(h^{-1}(w))$ is a $\vartheta_0$-discrete subset of $Y$. Then $f^{-1}(w) = \nu\cdot B$ is $\vartheta_0$-discrete in $L\subseteq [1,\gl]$, hence in $[1,\gl]$.

%\medskip
Finally, we need to show

\medskip
\textsc{Claim 5:} $f$ is algebra preserving.

We show that $A:=f^{-1}(w)\in\HH_n([1,\gl])$, for all $w\in T$. Again, we consider two cases.

\medskip
\hspace{\parindent}\textsc{Case 1:} $w\in S$. Then $A\subseteq M$ and $B:=g^{-1}(w) \in \HH_{n}([1,\nu])$ by the induction hypothesis and preservation under the sum. Hence, $A\cap I_\ga=\nu\ga+ B\in \HH_n(I_\ga)$, for each $\ga<\gw^\mu$.
We need to show that $A$ is u.p.\ in $\gb$ for each $\gb\leq\gl$.

If $\gb\in M$ then $\gb\in I_\ga$, for some $\ga<\gw^\mu$, and  $\gb=\nu\ga+\gb'$, for some $\gb'\in [1,\nu]$. $B$ is u.p.\ in $\gb'$, hence $A\cap I_\ga=\nu\ga+B$ is u.p.\ in $\gb$. Therefore $A$ is u.p.\ in $\gb$.

If $\gb\notin M$, then $\gb\in L$, that is, $\gb=\nu\ga$, for some $\ga\in\Lim$. For some $\gd>0$ there holds $\ga=\gw\gd$. If $\gd$ is a successor than $\gb = \nu(\gy+\gw)=\nu\gy+\nu\gw$, for some $\gy$. Then $A\cap \gb$ can be represented as a transfinite word $CB^\gw$, for some $C$. Hence, $A$ is u.p.\ in $\gb$.

Suppose $\gd$ is a limit ordinal. Then $A\cap \nu\gw$ can be represented as the word $0B^\gw$. Since $A\cap L=\emptyset$ we obtain that $A\cap \gb$ is  $(0B^\gw)^\gd$, which periodic.

\medskip
\hspace{\parindent}\textsc{Case 2:} $w\in Q$. By the induction hypothesis $h^{-1}(w)\in \HH_{n-1}([1,\mu])$. Hence, by definition, $B:=\ell^{-1}(h^{-1}(w))\in \HH_n([1,\gw^\mu])$. It follows that $\nu\cdot B\in \HH_n([1,\gl])$.
\ep

\section{Completeness}

 Now we are ready to prove our main theorem.

\bt \label{topcomp}
\benr
\item For each $n<\gw$ and $\Omega\geq \gw_{n+1}$, $$\Log(\Omega;\vartheta_0,\dots,\vartheta_n;\HH_n(\Omega))=\GLP_n;$$
\item  For $\Omega\geq\ge_0$,  $\Log(\Omega;(\vartheta_n)_{n<\gw};\HH_n(\Omega))=\GLP$.
\eenr
\et

\bp\  Suppose $\GLP_n\nvdash\varphi$. Obviously,  $\J_n\nvdash M^+(\varphi)\to\varphi$. Then there exists a finite $J_n$-tree $\cT$ such that $\cT\nmodels M^+(\varphi)\to\varphi$. By Lemma~\ref{lkey} there is an ordinal $\gl<\gw_{n+1}$ and a $J_n$-morphism $f:\cX\tto \cT$, where $\cX=([1,\gl],\vartheta_0,\dots,\vartheta_n;\HH_n([1,\gl]))$. By Theorem~\ref{t:main-1} we have $\cX\nmodels\varphi$. Since $\cX$ is isomorphic to a subframe of $(\Omega;\vartheta_0,\dots,\vartheta_n;\HH_n(\Omega))$ we have the required.
\ep

%\bibliographystyle{plain}
%\bibliography{ref-all2}

\section{Appendix: Proof of Theorem~\ref{t:main-1}}

Let $\tilde\D(A)$ abbreviate $X\setminus\D(X\setminus A)$.
Obviously, $x\in\tilde\D(A)$ iff $A$ contains some punctured
neighborhood of $x$.
\bl \label{l:j34} Assume $f:\cX\to \cT$ is a $J_n$-morphism. Then, for any hereditary $(k+1)$-root $w\in T$, \[
f^{-1}(R_k^*(w)\cup\{w\}) \subseteq \tilde\D_k(f^{-1}(R_k^*(w))).
\leqno (*) \]
\el

\bp\  Assume $a\in
f^{-1}(R_k^*(w)\cup\{w\})$. We have to construct a punctured
neighborhood of $a$ contained in $f^{-1}(R_k^*(w))$. Consider
$$U:=f^{-1}(R_k^*(w)\cup\{w\})=f^{-1}(R_k^*(w))\cup f^{-1}(\{w\}).$$
By the first part of $(j_3)$, $U$ is a neighborhood of $a$. If $a\in
f^{-1}(R_k^*(w))$ then $V:=f^{-1}(R_k^*(w))$ is a neighborhood of
$a$ by the second part of $(j_3)$, so $V-\{a\}$ is as required. If
$a\in f^{-1}(w)$ then by $(j_4)$ there is a neighborhood $V_a$ such
that $V_a\cap f^{-1}(w)=\{a\}$. Then, $$V_a\cap U= (V_a\cap
f^{-1}(R_k^*(w)))\cup \{a\}$$ is a neighborhood of $a$. Then,
$(V_a\cap U)\setminus\{a\}$ is a punctured neighborhood of $a$
contained in $f^{-1}(R_k^*(w))$.
\ep

\paragraph{Proof of Theorem~\ref{t:main-1}.}
Let $\cX$ be an ambience, $\cT$ a ${J}_n$-tree, $f:\cX\tto \cT$ an algebra preserving surjective
${J}_n$-morphism. We show that, for all $\mathcal{L}_n$-formulas  $\varphi$,
$\cX\models\varphi$ implies $\cT\models M^+(\varphi)\to\varphi$.

Suppose $\cT\nmodels M^+(\varphi)\to\varphi$. Then for some valuation
$\nu$ on $\cT$ and some point $w\in T$ (assume without loss of
generality that $w$ is the root of $\cT$) we have that
$w\in\nu(M^+(\varphi))$ but $w\not\in\nu(\varphi)$. Consider a
valuation $\nu'$ on $\cX$ by taking $\nu'(p):=f^{-1}(\nu(p))$.

\bl
For all subformulas $\theta$ of $\varphi$, we have
$\nu'(\theta)=f^{-1}(\nu(\theta))$.
\el

\bp\
We argue by induction on the complexity of $\theta$. If $\theta$ is
a propositional letter, the claim follows from the definition of
$\nu'$. The case of propositional connectives is easy.

If $\theta=[n]\psi$, then the claim follows by condition $(j_1)$ of $f$ being a $J_n$-morphism.

Suppose $\theta=[k]\psi$ for some $k<n$. To show that
$\nu'(\theta)\subseteq f^{-1}(\nu(\theta))$ assume
$x\in\nu'(\theta)$. Then there exists a $U\subseteq X$ such that
$\{x\}\cup U\in\tau_k$ and $U\subseteq\nu'(\psi)$. By IH we obtain
$U\subseteq f^{-1}(\nu(\psi))$. Hence $f(U)\subseteq
f(f^{-1}(\nu(\psi)))=\nu(\psi)$. By $(j_2)$, the set $f(\{x\}\cup
U)=\{f(x)\}\cup f(U)$ is an $R_k$-upset and so $R_k(f(x))\subseteq
f(U)\subseteq\nu(\psi)$. It follows that $f(x)\in\nu([k]\psi)$. In
other words, $x\in f^{-1}(\nu(\theta))$.

For the converse inclusion suppose $x\in f^{-1}(\nu(\theta))$, that
is,  $f(x)\models [k]\psi$. We must show $x\in \nu'(\theta)$. By the
induction hypothesis,
$$\nu'(\theta)=\tilde\D_k(\nu'(\psi))=\tilde\D_k(f^{-1}(\nu(\psi))).$$

Let $v\in T$ be a hereditary $(k+1)$-root such that $v=f(x)$ or $v
R_{k+1} f(x)$. Since $\cT$ is a $J_n$-tree,
$R_k(v)=R_k(f(x))$. Thus $v\models[k]\psi$. We also have $v\models
M^+(\varphi)$. In particular, $v\models [k]\psi\to[k']\psi$ for any
$k'$ with $k\leq k'\leq n$ and hence $v\models[k']\psi$. It follows
that for each $k'$ between $k$ and $n$ we have
$R_{k'}(v)\subseteq\nu(\psi)$. Therefore $R^*_{k}(v)\subseteq
\nu(\psi)$ and hence $f^{-1}(R^*_{k}(v))\subseteq
f^{-1}(\nu(\psi))$. By the construction of $v$, $x\in
f^{-1}(R^*_{k}(v)\cup\{v\})$. Hence,  by Lemma \ref{l:j34},
$$x\in \tilde\D_k(f^{-1}(R^*_{k}(v)))\subseteq\tilde\D_k(f^{-1}(\nu(\psi))),$$
as required. \ep

From this lemma we obtain $y\notin
\nu'(\varphi)=f^{-1}(\nu(\varphi))$, for any $y$ with $f(y)=w$.
Consequently, $\cX\nmodels\varphi$.
\ep

\ignore{For the (if) part of the theorem we assume $\cX\nmodels \phi$. Then $\cX\models\not\phi$ and by the (only if) part $\cT\models M^+(\neg\phi)\to \neg\phi$. If $\cT\models M^+(\phi)\to\phi$ then we conclude $\cT\nmodels M^+(\phi)\land M^+(\neg\phi)$. Hence, there is a valuation $\nu$ on $\cT$ and a formula $\psi$  such that $\nu([i]\psi))\nsubseteq \nu([j]\psi))$, for some $i<j\leq n$. Then $B:=f^{-1}(\nu(\psi))$ is admissible and $d_i(B)\nsubseteq d_j(B)$.}

\end{document}